\input amstex
\documentstyle{amsppt}
\magnification 1200
\hsize=14truecm
\vsize=23truecm
\hfuzz 1pt
\tolerance 1000

\def\C{{\Bbb C}}
\def\Z{{\Bbb Z}}
\def\Q{{\Bbb Q}}
\def\Aa{{\Bbb{A}}}
\def\Pp{{\Bbb{P}}}
\def\Ss{{\Bbb{S}}}
\def\cQ{{\Bbb Q}}
\def\Pr{{\bold P}}
\def\A{{\Cal{A}}}
\def\uC{{\Cal C}}

\def\H{{{\Cal H}}}

\def\L{{{\Cal L}}}

\def\O{{\Cal O}}
\def\PP{{\Cal P}}
\def\P{{\Cal P}}

\def\cS{{\Cal S}}

\def\ZZ{{\Cal Z}}
\def\sO{{\Cal O}}
\def\sD{{\Cal D}}
\def\sM{{\Cal M}}
\def\sN{{\Cal N}}
\def\sE{{\Cal E}}

\def\W{{{\Lambda }}}
\def\M{{\overline{\Cal M}}}

\def\Res{{\operatorname {Res}}}
\def\Sym{{\operatorname {Sym}}}
\def\Proj{{\operatorname {Proj}}}
\def\weight{{\operatorname {weight}}}
\def\Aut{{\operatorname {Aut}}}
\def\Spec{{\operatorname {Spec}}}
\def\cl{{\operatorname {cl}}}

\def\Tr{\operatorname{Tr}}
\def\Div{\operatorname{Div}}

\def\Supp{\operatorname{Supp}}
\def\Sp{\operatorname{Spec}}
\def\Hom{\operatorname{\underline{Hom}}}

\def\LL{{\L\negthickspace\L}}
\def\CP{{\C\bold P}}
\def\Discr{{\Sigma}}

\def\ls{\leqslant}
\def\gs{\geqslant}

\topmatter
\title{Hurwitz numbers and intersections on moduli spaces of curves}
\endtitle

\author Torsten Ekedahl$^*$,
Sergei Lando$^{\dag,}$\footnote
{Partly supported by the RFBR grant 01-98-00555 and by the KVA grant
for Cooperation between Sweden and the former Soviet Union.\kern 5cm}, 
Michael Shapiro$^\ddag$, and 
Alek Vainshtein$^\natural$
\endauthor
\affil
$^*$ Dept. of Math., University of Stockholm,
S-10691, Stockholm, {\tt teke\@matematik.su.se}\\
$^\dag$ Higher College of Math.,
Independent University of Moscow, and Institute for System Research RAS, 
%partly supported by the RFBR grant 01-98-00555, 
{\tt lando\@mccme.ru}\\
$^\ddag$ Department of Mathematics,
Royal Institute of Technology, S-10044, Stockholm, 
{\tt mshapiro\@math.kth.se}\\
$^\natural$ Dept. of Math.
and Dept. of Computer Science, University of Haifa, Haifa 31905,
{\tt alek\@mathcs.haifa.ac.il}
\endaffil
\leftheadtext{ Ekedahl, Lando, Shapiro, Vainshtein}
\rightheadtext{Hurwitz numbers and moduli spaces of curves}
\endtopmatter

\document
\heading 1. Introduction \endheading

\subheading{1.1. Topological classification of
ramified coverings of the sphere}

\noindent For a compact connected genus $g$ complex curve $C$
let $f\:C\to\CP^1$ be a meromorphic function.
We treat this function as a ramified covering
of the sphere. Two ramified coverings $(C_1;f_1)$, $(C_2;f_2)$
are called {\it topologically equivalent\/} if there exists a
homeomorphism $h\:C_1\to C_2$ making the following diagram
commutative:
$$
\matrix
C_1&@>h>>&C_2\\
f_1\searrow && \swarrow f_2\\
&\CP^1\endmatrix
$$
The critical values of topologically equivalent functions,
i.e., the ramification points of the coverings,
coincide, as do the genera of the covering curves.
In his famous paper~\cite{H}  Hurwitz initiated the
topological classification of such coverings in the case when
exactly one of the ramification points is degenerate, and the
remaining points are nondegenerate.
Below we refer
to the degenerate ramification point as ``infinity'',
and its preimages are called ``poles''.

The space of all
biholomorphic equivalence classes of
generic (that is, with distinct critical values)
meromorphic functions on genus $g$ curves
with the prescribed orders $k_1,\dots,k_n$
of poles carries a natural topology, a natural complex structure,
and is connected (see~\cite{N}). The  $\C$-dimension of this
 space is $\mu=K+n+2g-2$, where $K=k_1+\dots+k_n$
is the number of folds in the covering.
The additive group $\C$ of complex numbers
acts on this space by adding a constant $c\in\C$ to
all meromorphic functions. The {\it Hurwitz space\/}
is the space of orbits with respect to this action.
There is a unique meromorphic function in each orbit
such that the sum of its finite critical values is zero.
Hence, the Hurwitz space is naturally
identified with the space of meromorphic functions
with zero sum of finite critical values.
We denote the latter space by $H_{g;k_1,\dots,k_n}$.

For a given set of orders $k_1,\dots,k_n$,
the number of the equivalence classes of
topologically nonequivalent ramified coverings
with these orders of poles and prescribed nondegenerate
ramification points is finite; it is
called the {\it Hurwitz number\/} and denoted $h_{g;k_1,\dots,k_n}$.
Hurwitz posed the problem of finding $h_{g;k_1,\dots,k_n}$ explicitly.
Below we express the Hurwitz number
in terms of intersection numbers for the Chern classes
of certain bundles on the moduli space of
complex curves with $n$ marked points.

Let $\M_{g;n}$ denote the Deligne--Mumford compactification
of the moduli space of complex curves with $n$ marked points,
and let $\L_i$ be the line bundle on $\M_{g;n}$
whose fiber at a point $(C;x_1,\dots,x_n)\in\M_{g;n}$
coincides with the cotangent space to $C$ at $x_i$.
The (first) Chern class of
such a bundle is denoted by $\psi_i=c_1(\L_i)$.

%More exactly, $\psi_i$ is the first Chern class of
%$\sigma^*_i(\omega_{\pi_{n+1}})$, where $\pi_{n+1}:\ \M_{g;n+1}\to\M_{g;n}$
%is the morphism obtained by forgetting the $(n+1)$st marked point,
%$\omega_{\pi_{n+1}}$ is the relative dualizing sheaf, and $\sigma_1,\dots,
%\sigma_n$ are the natural sections of $\pi_{n+1}$ (the image of
%$(C;z_1,\dots,z_n)$ under $\sigma_i$ is obtained by attaching a 3-pointed
%rational curve at $z_i$ and considering the remaining 2 points on that curve
%as $z_i$ and $z_{n+1}$).

The main goal of the present paper is the proof of the following
theorem.

\proclaim{Theorem 1.1}
For  $g=0$, $n\gs3$, or for $g\gs1$,
the Hurwitz number $h_{g;k_1,\dots,k_n}$ equals
$$
h_{g;k_1,\dots,k_n}=
\frac{(K+n+2g-2)!}{\#\Aut(k_1,\dots,k_n)}
\prod_{i=1}^n\frac{k_i^{k_i}}{k_i!}
\int_{\M_{g;n}}\frac{c(\W^\vee_{g;n})}{(1-k_1\psi_1)\dots(1-k_n\psi_n)}.
\tag1
$$
\endproclaim

Here $\#\Aut(k_1,\dots,k_n)$ is the number of automorphisms
of the $n$-tuple $(k_1,\dots,k_n)$,
$\W_{g;n}$ denotes
the Hodge  bundle of holomorphic 1-forms over $\M_{g;n}$
(a precise definition of this bundle is given below),
and $c(\W^\vee_{g;n})$ is the total Chern class of the dual bundle.
We admit as usual that the integral of a class
with degree different from the dimension of the
variety is zero.

Observe that the integral in~(1)
in case $g=0$, $n\gs3$ is equal to $K^{n-3}$
(see Sec.~2.1). It follows that the two exceptional cases $g=0$, $n=1,2$
are also covered by~(1),
 provided the integral is understood formally as
$K^{n-3}$.

Our main tool in the proof of the main theorem
is the Lyashko--Looijenga mapping described below.

\subheading{1.2. Lyashko--Looijenga mapping}
The Lyashko--Looijenga mapping (the LL mapping in the
text and $\LL$ in the formulas below) associates to
a $\C$-valued holomorphic function the unordered set of its
critical values (taking multiplicities into account).
It is a classical tool in the study of the geometry of
moduli spaces of meromorphic functions.
Assume we have a family of functions such that for
each member of the family the set of its critical
values is finite and contains the same number
of elements. Then the LL mapping can be considered
as a mapping from the family to the space of
monic (i.e. with leading coefficient~$1$)
polynomials in one variable. To do this  we associate to an
unordered set of complex numbers the monic polynomial
whose roots coincide with the elements of the set.
The most interesting situation occurs when
the LL mapping is of an algebraic nature, and one can compute
its multiplicity. We are going to show that this is
the case for the Hurwitz spaces.

The Riemann--Hurwitz formula implies that a
function belonging to the Hurwitz space
has $\mu=K+n+2g-2$ finite ramification points.
Recall that the sum of the finite critical values equals zero.
Therefore the LL mapping can be treated as the mapping to
$\C^{\mu-1}$,
where the target space is the
space of polynomials of the form
$$
t^{\mu}+d_2t^{\mu-2}+\dots+d_{\mu}, \tag2
$$
and the coefficients $d_i$ of these polynomials
form a natural set of coordinates in this space.
We denote the Lyashko--Looijenga mapping by $\LL$,
$$
\LL\:H_{g;k_1,\dots,k_n}\to\C^{\mu-1}.
$$

The multiplicity of the LL mapping
is closely related
to the enumeration problem for topological
types of ramified coverings.

\proclaim{Lemma 1.2}
$$
h_{g;k_1,\dots,k_n}= \frac{\deg\ \LL}{\#\Aut(k_1,\dots,k_n)}.
$$
\endproclaim

\demo{Proof}
Indeed, consider a  function $f\:C\to\CP^1$
with a prescribed set of $\mu$ distinct
critical values. Such function
determines a complex structure on~$C$.
A topological equivalence preserves the set of poles
of the function, and it can permute only poles
of the same order. If the poles are fixed, then,
according to the Riemann theorem,
the meromorphic function is unique. Thus, the
multiplicity of the LL mapping is simply
$\#\Aut(k_1,\dots,k_n)$ times the number of topological
types of meromorphic functions.
\qed
\enddemo

\subheading{1.3. Outline of the proof}
In the present section we describe briefly the most important steps
in the proof of the main theorem.

We start with the definition of the space $\PP$ of generalized principal parts.
This space is considered as a cone (in the sense of~\cite{Fu})
over the moduli space $\M_{g;n}$ of stable curves
with $n$ marked points. The fiber of $\PP$ at a point $(C;x_1,\dots,x_n)$
consists essentially of $n$-tuples of principal parts
of meromorphic germs with poles of order $k_i$
at the marked points. This cone is simply a vector bundle
if the orders $k_i$ of all poles are equal to~1, but
it looses the linearity in the case
of larger orders.

Associating to a meromorphic function the $n$-tuple
of its principal parts at the marked points
we define an embedding of the Hurwitz space to $\PP$.
If $g=0$, then the completed Hurwitz space simply coincides
with the space of generalized principal parts. For $g>0$,
the Hurwitz space is a subcone in the cone $\PP$.
The principal parts $(p_1,\dots,p_n)$ of a meromorphic function
$f\:(C;x_1,\dots,x_n)\to(\CP^1,\infty)$
at the marked points $x_i$ must satisfy the requirement
$$
\Res_{x_1} p_1\omega+\dots+\Res_{x_n} p_n\omega=0 \tag3
$$
for any holomorphic 1-form $\omega$ on $C$.
We {\it define\/} the completed Hurwitz space as the closure
in $\PP$ of the set of principal parts on {\it smooth curves\/}
satisfying requirement~(3). The completed Hurwitz
space has a natural structure of a cone over $\M_{g;n}$.

The LL mapping extends naturally to the completed Hurwitz
space. Its multiplicity can be expressed in terms of the top
Segre class of this cone. On the other hand, this top Segre
class can be expressed
in terms of the total Segre class of the cone of generalized principal parts
and of the total Chern class of the Hodge bundle over $\M_{g;n}$.
The main difficulty here arises from the fact that
relation~(3) is satisfied not only by the points of the
Hurwitz space, but as well by some other $n$-tuples of principal
parts coming from meromorphic functions on singular curves.
The {\it total\/} Segre class
takes into account these additional components.
Fortunately, however, their impact on the {\it top\/}
Segre class, which is the one we are interested in,
can be shown to be trivial, and we arrive at the integral~(1).

\subheading{1.4. Previous research}
The study of topological types of ramified
coverings of the sphere was initiated by A.~Hurwitz
in~\cite{H}. In particular, he suggested there,
without a complete proof, an explicit formula for the number
of genus $0$ ramified coverings of the sphere.
We reproduce this formula in Sec.~2.
In~1995 it was partially (in the case
of simple poles) rediscovered by physicists~\cite{CT}.
Inspired by their result  I.~Goulden and
D.~Jackson~\cite{GJ1} came to a proof of the general Hurwitz formula
(which they were not aware of). Their approach, mainly
purely combinatorial, is presumably close to
the original Hurwitz way of reasoning (see the reconstruction
of Hurwitz's proof in~\cite{S}). The further
development of the same approach resulted in
obtaining explicit enumeration of ``simple''
toric coverings (these formulas are
also reproduced in Sec.~2), and a number of other partial
results and interesting conjectures in this direction~\cite{SSV,GJVn}.

On the other hand, Mednykh~\cite{M1,M2} gave, in a sense,
a complete answer to the enumeration problem under
consideration. Unfortunately, the way of presentation
the results does not allow one to extract essential
information about the behavior of the numbers
and the underlying geometry remains covered.

The Lyashko--Looijenga mapping was involved in the
subject by V.~I.~Ar\-nold~\cite{A1}. The mapping itself
was introduced by Lyashko (unpublished, see~\cite{A2})
and, independently, by Looijenga~\cite{L} as a main tool
of investigating the topology of the complement
to the discriminant in the space of versal
deformations of simple singularities.
Arnold reinterpreted the results of Lyashko and Looijenga
for the singularities of the $A_\mu$ series as
the topological classification of generic rational ramified
coverings of the sphere with one pole and
extended them to the case of two poles. Further
exploitation of the same tools~\cite{GL} led
to a new proof of Hurwitz's formula in all
its generality. Our present approach is close
to that of~\cite{GL}, although we consider
both a different compactification of the moduli
space of curves, and a different fibration over this
space. We do not
know a direct generalization of the construction
of~\cite{GL} to higher genera.

On degenerate genus $0$ coverings see~\cite{GJ2,Z,LZ}.

The last years' flash of interest to the topology of
Hurwitz spaces is due to the fact that it
can be treated as one
of the easiest examples of quantum cohomology
calculations~\cite{KM,V1,GP}.
It is also worth mentioning that the LL mapping
is the main tool in Dubrovin's construction of
various Frobenius structures in Hurwitz spaces~\cite{D}.

The results of the present paper were announced in~\cite{ELSV}
(with a number of mistakes that are corrected here).
The main idea of the present proof follows that of~\cite{ELSV},
although technical details are different.
Since~\cite{ELSV} has been published, an independent proof
of the formula~(1) for the case of all simple poles,
$k_i=1$, appeared in~\cite{FnP}. More recently, a proof of (1) for the 
general case based on the same ideas of virtual localization as the proof 
presented in~\cite{FnP}, appeared in~\cite{GV}.

The authors are sincerely grateful to V.I.Arnold who introduced the last
three authors to this subject.
We are indebted to B.~Dubrovin, L.~Ernstr\"om, C.~Faber,
A.~Gorodentsev, A.~Levin, Yu.~Manin,
S.~Natanzon, R.~Pandharipande, M.~Rosellen and,
especially, B.~Shapiro and R.~Vakil
for many fruitful discussions.
During the work on this paper
the authors enjoyed the hospitality of the University
of Stockholm (S.L., A.V.), the University of Rennes (S.L.),
the University of Haifa (M.S.),
the Max-Planck-Institut f\"ur Mathematik (S.L), and the KTH (S.L., A.V.).

\heading {2. Calculations}
\endheading

\subheading{2.1. The genus zero case}
In the genus zero case the integral in the right hand side
of the main formula~(1) has the form
$$
\int_{\M_{0;n}}\frac1{(1-k_1\psi_1)\dots(1-k_n\psi_n)}, \tag4
$$
and in order to compute it it is sufficient to find the
intersection numbers for the monomials
$$
\langle\tau_{m_1}\dots\tau_{m_n}\rangle_0=
\int_{\M_{0;n}}\psi_1^{m_1}\dots \psi_n^{m_n}, \tag5
$$
which are special cases, for $g=0$, of the monomials
$$
\langle\tau_{m_1}\dots\tau_{m_n}\rangle_g=
\int_{\M_{g;n}}\psi_1^{m_1}\dots\psi_n^{m_n}.
$$
Intersection numbers~(5) are nonzero if and only if
the degree of the integrand coincides with the
dimension of the base,
$$
m_1+ m_2+\dots+ m_n=n-3.
$$
They are totally determined by the initial condition
$\langle\tau_0^3\rangle_0=1$ and the genus zero case of the
{\it string equation}
$$
\langle\tau_0\prod_{i=1}^n\tau_{m_i}\rangle_g
=\sum_{j=1}^n\langle\tau_{m_j-1}\prod_{i\ne j}\tau_{m_i}\rangle_g \tag6
$$
(Witten's theorem for genus~0, see~\cite{K,W}).
An easy computation gives
$$
\langle\tau_{m_1}\dots\tau_{m_n}\rangle_0=
\frac{(n-3)!}{m_1!\dots m_n!}.
$$
This formula was known to physicists for about
the last two decades.

Hence, the coefficient at $k_1^{m_1}\dots k_n^{m_n}$
in the expansion of~(4) equals zero if
$m_1+\dots+m_n\ne n-3$,
and equals
$$
\frac{(n-3)!}{m_1!\dots m_n!}
$$
otherwise. Therefore, (4) is nothing but
the expansion of $(k_1+\dots+k_n)^{n-3}$.
We thus arrive to the expression
$$
h_{0;k_1,\dots,k_n}=\frac{(K+n-2)!}{\#\Aut(k_1,\dots,k_n)}
\prod_{i=1}^n\frac{k_i^{k_i}}{k_i!}\cdot K^{n-3}, \tag7
$$
where $K=k_1+\dots+k_n$,
which is known since the pioneering work of
Hurwitz~\cite{H}.

\subheading{2.2. Calculation of some intersection numbers
for higher genera}
For $g\gs1$ the factor $c(\W_{g;n}^\vee)$ in the integrand in~(1)
is not 1 any more. It has the form
$$
c(\W_{g;n}^\vee)=1-\lambda_1+\dots+(-1)^g\lambda_g,
$$
where $\lambda_i\in H^{2i}(\M_{g;n})$.
Hence, the calculation of the integral is reduced
to the calculation of the integrals of monomials
of the form
$$
\int_{\M_{g;n}}\psi_1^{m_1}\dots \psi_n^{m_n}\lambda_i
$$
for
$$
m_1+\dots +m_n=3g+n-3-i
$$
(for other monomials the integral vanish). These integrals
are called {\it Hodge integrals} (see~\cite{FP1, FP2}),
and their values are known in some special cases.

In particular, it is known that

\proclaim{Theorem 2.1} {\rm \cite{FP2}}
$$
\int_{\M_{g;n}}\psi_1^{m_1}\dots \psi_n^{m_n}\lambda_g
= {2g+n-3\choose m_1\dots m_n}b_g,
$$
where $b_g$ is a constant independent of $m_j$ and equal to
$$
b_g=\cases
1, & g=0,\\
\int_{\M_{g;1}}\psi_1^{2g-2}\lambda_g, & g>0.
\endcases
$$
\endproclaim

From this result we immediately get the following corollary.

\proclaim{Corollary 2.2}
The term of the lowest degree in  $k_1,\dots,k_n$
in the integral~{\rm(1)} equals, up to a nonzero
constant factor,
$$
(k_1+\dots+k_n)^{2g+n-3}.
$$
The factor is $\int_{\M_{g;1}}\psi_1^{2g-2}\lambda_g$.
\endproclaim

In particular, for $g=1$ this statement
entirely determines the ``$\lambda$-part''
of the integral in the right-hand side of~(1).
The ``$\lambda$-free'' part of the integral
is, in principle, described by Witten's conjecture
(Kontsevich's theorem,~\cite{W,K}) for all genera.
In addition to the string equation~(6), the description includes
the {\it dilaton equation}
$$
\langle\tau_1\prod_{i=1}^n\tau_{m_i}\rangle_g=
(2g-2+n)\langle\prod_{i=1}^n\tau_{m_i}\rangle_g. \tag8
$$
Together with the initial value
$$
\langle\tau_1\rangle_1=\int_{\M_{1;1}}\lambda_1=\frac1{24},
$$
this gives a complete description of genus~one Hurwitz numbers.
In particular, it proves the following elegant closed formula
 conjectured in \cite{GJVn} and first proved in \cite{GJ3, V2}.

\proclaim{Theorem 2.3}
For $g=1$
$$
h_{1;k_1,\dots,k_n}=
\frac{(K+n)!}{24\#\Aut(k_1,\dots,k_n)}
\prod_{i=1}^n\frac{k_i^{k_i}}{k_i!}
\bigg(K^n-\sum_{i=2}^n(i-2)!e_iK^{n-i}-K^{n-1}\bigg),
$$
where $e_i=e_i(k_1,\dots,k_n)$ is the $i$th elementary
symmetric polynomial in $k_1,\dots,k_n$. 
%$e_1=K=k_1+\dots+k_n$.
\endproclaim

\demo{Proof} It follows from the main theorem and
Corollary~2.2 that it is enough to prove that
$\langle\tau_{m_1}\dots\tau_{m_n}\rangle_1$ equals the coefficient
at the monomial $k^{m_1}\dots k^{m_n}$ in the expression
$$
\frac1{24}\left(K^n-\sum_{i=2}^n(i-2)!e_iK^{n-i}\right).\tag9
$$
The latter coefficient is equal to 
$$
\alpha_{m_1,\dots,m_n}=
\frac1{24M}\left(n!-\sum_{i=2}^n(i-2)!(n-i)!e_i(m_1,\dots,m_n)\right),
$$ 
where $M=m_1!\dots m_n!$.

First, assume that all $m_j$'s are distinct from zero, and   hence
$m_1=\dots=m_n=1$. Then, by the dilaton equation,
$$
\langle\tau_{1}\dots\tau_{1}\rangle_1=\frac{(n-1)!}{24},
$$
while the corresponding coefficient in~(9) equals
$$
\alpha_{1,\dots,1}=
\frac1{24}\left(n!-\sum_{i=2}^n (i-2)!(n-i)!{n\choose i}\right)=
\frac{(n-1)!}{24}.
$$

Now, assume that at least one  among $m_j$'s is equal to zero, say,
$m_1=0$. Define $J=\{j\: m_j\gs1\}$.
It is easy to see that
$$\multline
e_i(m_2,\dots,m_j-1,\dots,m_n)\\
=e_i(m_2,\dots,m_j,\dots,m_n)-
e_{i-1}(m_2,\dots,0,\dots,m_n)
\endmultline
$$
for $j\in J$ and
$$
\sum_{j=2}^nm_je_{i-1}(m_2,\dots,0,\dots,m_n)=ie_i(m_2,\dots,m_j,\dots,m_n).
$$
Therefore  we get
$$\multline
\sum_{j\in J}\alpha_{m_2,\dots,m_j-1,\dots,m_n}\\
\shoveleft{=\frac1{24M}\sum_{j\in J}m_j
\bigg((n-1)!-\sum_{i=2}^{n-1}(i-2)!(n-1-i)!}\\
\shoveright{\times\left(e_i(m_2,\dots,m_j,\dots,m_n)-
e_{i-1}(m_2,\dots,0,\dots,m_n)\right)\bigg)}\\
=\frac1{24M}\bigg(n!-\sum_{i=2}^{n-1}(i-2)!(n-i)!
e_i(m_2,\dots,m_j,\dots,m_n)\bigg)
= \alpha_{0,m_2,\dots,m_n},
\endmultline
$$
since $e_i(m_2,\dots,m_j,\dots,m_n)=e_i(0,m_2,\dots,m_n)$ and
$e_n(0,m_2,\dots,m_n)=0$.

So, coefficients $\alpha_{m_1,\dots,m_n}$ obey the same string equation as
$\langle\tau_{m_1}\dots\tau_{m_n}\rangle_1$ and coincide with the
latter for initial values of $m_j$'s. Therefore, they coincide for
arbitrary values of $m_j$'s.
\qed
\enddemo

Another application of the main theorem is the following generating
function for Hodge integrals over moduli spaces with one
marked point found first in~\cite{FP1}.

\proclaim{Theorem 2.4}
$$
1+\sum_{g=1}^\infty
t^{2g}\int_{\M_{g;1}}
\frac{k^g-k^{g-1}\lambda_1+\dots+(-1)^g\lambda_g}{1-k\psi}=
\left(\frac{t/2}{\sinh\ t/2}\right)^{k+1}.
$$
\endproclaim

\demo{Proof} Follows immediately from the main theorem and the
generating
function for the Hurwitz numbers $h_{g;k}$ obtained in \cite{SSV}.
\qed
\enddemo

Finally, as an immediate corollary of the main theorem
we obtain the following statement, which generalizes
a conjecture from~\cite{GJVn}.

\proclaim{Theorem 2.5}
The factor in $h_{g;k_1,\dots,k_n}$
given by the integral is the sum of
homogeneous symmetric polynomials in $k_1,\dots,k_n$
of degrees $n+3g-3,n+3g-4,\dots,n+2g-3$.
\endproclaim

For other applications of the main theorem see~\cite{GJVk}.

\heading {3. Completed Hurwitz spaces}
\endheading

\subheading{3.1. Cones and projective cones}
Let $\Ss^\cdot=\Ss^0\oplus \Ss^{1}\oplus \Ss^{2}\oplus\dots$
be a graded sheaf of $\O_X$-algebras over a scheme $X$.
%$q_i\in\Q$, $0<q_1<q_2<\dots$.
We suppose that $\Ss^\cdot$
is locally finitely generated and
the sheaf of fields of quotients  of  $\Ss^\cdot$
coincides with that of $\Ss^1$.
Note that this requirement is weaker than the one used in~\cite{Fu}
that $\Ss^\cdot$ itself is locally finitely generated
as an $\O_X$-algebra by $\Ss^1$.

Given $\Ss^\cdot$, we define two schemes over $X$,
the {\it cone\/}
$
\cS=\Spec(\Ss^\cdot)
$
endowed with the action of the group $\C^*$
due to the grading in $\Ss^\cdot$,
and the {\it weighted projective cone\/}
$
P\cS={\Proj}(\cS)
$
(the {\it projectivization\/} of $\cS$).
Points of $P\cS$ correspond to non-trivial orbits
of the $\C^*$-action on $\cS$.
The projective cone $P\cS$ is endowed with the canonical
line bundle $\O(1)$ of  quasihomogeneous
rational functions of weight~1.

Any vector bundle $E$ has a natural structure of a cone determined
by its structure sheaf $\Sym(E^\vee)$
endowed with the natural (integer-valued)
grading by degrees.

For two cones $\cS_1,\cS_2$ over $X$ defined by the sheaves 
$\Ss^\cdot_1,\Ss^\cdot_2$
their {\it direct sum\/} $\cS_1\oplus \cS_2$ is defined by the graded sheaf
$\Ss^\cdot_1\otimes\Ss^\cdot_2$.

\subheading{3.2. Moduli space of stable curves}
The base space of most of the bundles and cones considered in the
present paper will be the {\it moduli space\/} $\M_{g;n}$
of stable curves with $n$ marked points. Its elements are
the biholomorphic equivalence
classes of connected compact curves of arithmetic genus $g$
with $n$ marked nonsingular points that are either smooth
or have at most ordinary double points as singularities
and that admit no continuous group of automorphisms.
The last assumption (the {\it stability\/}) is equivalent
to the requirement that each rational (of genus~0) irreducible component
of the curve contains at least~3 special (marked or double)
points, and each elliptic (of genus~1) irreducible component
contains at least one special point.

The moduli space $\M_{0;n}$, $n\gs3$ of rational curves is a smooth
and compact variety, while for $g\gs1$ $\M_{g;n}$ is a compact
orbifold. By the Deligne--Mumford theorem,
the space $\M_{g;n}$ is an irreducible projective variety of pure
dimension $\dim\ \M_{g;n}=n+3g-3$.
All the facts about moduli spaces we make use of can be found in~\cite{HMo}.

A meromorphic function on a singular curve
is a tuple of meromorphic functions, one for each
irreducible component of the curve,
such that the values of these functions on each two branches
meeting at a double point coincide at this point.

\subheading{3.3. The space of generalized principal parts at a point}
Fix a positive integer~$k$.
Two germs of meromorphic functions $f_1,f_2\:(\C,0)\to(\CP^1,\infty)$
with poles of order $k$ have the same principal part
if their difference $f_1-f_2$ has no pole at~0.
A {\it principal part\/} is an equivalence class of germs
with respect to this equivalence relation.
The set of all principal parts with poles of order
precisely $k$ carries a natural complex structure.
Below, we present a coordinate description of the
space of principal parts.

Let $x$ be the germ of a coordinate at the origin $0\in\C$.
Then a principal part can be written in the form
$$
\left(\frac{u}x\right)^k+
a_1\left(\frac{u}x\right)^{k-1}+\dots
+a_{k-1}\frac{u}x,\qquad u\ne0. \tag10
$$
This presentation is not unique: it depends on the choice of
the parameter $u$,
and there are $k$ possibilities to make this choice.
Hence, the space of expressions of the form~(10)
covers the space of principal parts with multiplicity $k$.
The group of the covering is
$\Z/k\Z$, and it acts on the space of expressions~(10)
according to the rule
$$
(u,a_1,a_2,\dots,a_{k-1})\mapsto
(\zeta u,\zeta a_1,\zeta^2 a_2,\dots,\zeta^{k-1}a_{k-1}), \tag11
$$
where $\zeta$ is a primitive root of unity of degree~$k$
generating $\Z/k\Z$ as a subgroup in~$\C^*$.
We call the parameters $u,a_1,\dots,a_{k-1}$ the
{\it twisted Laurent coefficients of order~$k$}.

The group $\C^*$ of nonzero complex numbers acts on the space of
principal parts by multiplication.
Up to the above action of $\Z/k\Z$, this action of $\C^*$
can be written in coordinates as
$$\multline
c\:\left(\frac{u}x\right)^k+
a_1\left(\frac{u}x\right)^{k-1}+\dots+a_{k-1}\frac{u}x\\
\mapsto
c\left(\left(\frac{u}x\right)^k+
a_1\left(\frac{u}x\right)^{k-1}+\dots+
a_{k-1}\frac{u}x\right)\\
=\left(\frac{\eta u}x\right)^k+
(\eta a_1)\left(\frac{\eta u}x\right)^{k-1}+\dots+
(\eta^{k-1}a_{k-1})\frac{\eta u}x,
\endmultline
$$
where $\eta=c^{1/k}$.
Therefore, it induces the following choice
of {\it weights} of the parameters $u,a_j$:
$$
\weight(u) =\frac1k,\qquad
\weight(a_j)=\frac{j}k. \tag12
$$
Consider the algebra of polynomials in
twisted Laurent coefficients
invariant with respect to the action~(11) of $\Z/k\Z$.
The weights of the coordinates endow this algebra
with an integer-valued grading.
Denote the spectrum of this algebra by $P^k$.

\proclaim{Lemma 3.1}
{\rm 1)}
The graded algebra of $\Z/k\Z$-invariant polynomials is
independent of the choice of the local coordinate~$x$.

{\rm 2)}
The space of principal parts with poles of order $k$
is naturally embedded in $P^k$ as an open dense subset.

{\rm 3)}
Its complement $A^{k-1}\subset P^k$ carries a
natural structure of the quotient of
the complex vector space $\widetilde A^{k-1}\cong\C^{k-1}$
modulo the $\Z/k\Z$-action.

{\rm 4)}
The twisted Laurent coefficients $a_1,\dots,a_{k-1}$
are linear coordinates on $\widetilde A^{k-1}$.
The action of $\Z/k\Z$ is diagonal in these coordinates.
\endproclaim

\demo{Proof} A direct calculation shows that a coordinate change
$x=\alpha_1\tilde x+\alpha_2\tilde x^2+\dots$ with
$\alpha_1\ne0$ causes a polynomial change of the variables
$u,a_j$ of the form
$$
\align
\tilde u &=\frac{u}{\alpha_1^{k}},\\
\tilde a_1&=a_1+\gamma_{11}u,\\
\dots&\hdots\hdots\hdots\\
\tilde a_{k-1}&=a_{k-1}+\gamma_{k-1,1}a_{k-2}u+\dots+
\gamma_{k-1,k-1}u^{k-1}
\endalign
$$
for some constants $\gamma_{ij}$.
Evidently, this action  preserves the $\Z/k\Z$-invariance
and the grading.
Hence, the algebra of invariant polynomials is well-defined.

A principal part is given by an element in $P^k$
with $u\ne0$, and the set $u\ne0$ is open and dense in $P^k$.
Its complement $A^{k-1}\subset P^k$ is the spectrum
of the quotient algebra of $\Z/k\Z$-invariant polynomials
modulo the ideal of $u$-divisible polynomials.

Expression~(10) is linear in the coefficients $a_j$.
After setting $u=0$
the vector space structure on the space
$\widetilde A^{k-1}$ of these coefficients
is introduced in the obvious way
$$
c'(a'_1,\dots,a'_{k-1})+
c''(a''_1,\dots,a''_{k-1})=
(c'a'_1+c''a''_1,\dots,c'a'_{k-1}+c''a''_{k-1}).
$$
Its independence of the choice of the coordinate $x$
is the result of a direct computation.
The space $A^{k-1}$ is the quotient of $\widetilde A^{k-1}$
modulo the action~(11) restricted to~$\widetilde A^{k-1}$.

In order to prove that the coordinates in $A^{k-1}$
can be chosen invariantly,
fix a coordinate $x$ and consider a smooth holomorphic
1-parameter family $\gamma$ of principal parts through $a$,
$\gamma:(\C^1,0)\to (P^k,a)$, such that $\gamma(\tau)\notin A^{k-1}$
for $\tau\in\C^1$, $\tau\ne0$. Then the elements
of the family can be written in the form
$$
\left(\frac{u(\tau)}x\right)^k+
a_1(\tau)\left(\frac{u(\tau)}x\right)^{k-1}+\dots
+a_{k-1}(\tau)\frac{u(\tau)}x,\qquad u(\tau)\to0\ \text{as $\tau\to0$,}
$$
with $u(\tau),a_j(\tau)$ chosen uniquely up to the action~(11)
of $\Z/k\Z$.
One verifies immediately that the limit values
$a_1(0),\dots,a_{k-1}(0)$ depend neither on the choice of the
coordinate~$x$, nor on the choice of the 1-parameter family
through $a$.
\qed
\enddemo

Below, the term {\it generalized principal part of order~$k$}
means simply a point in $P^k$. Elements of $A^{k-1}$
are just generalized principal parts
{\it with zero leading coefficient}.

Now we are going to associate to a point in $P^k$
the $k$th tensor power of a tangent vector at $0\in\C^1$.
Let $L$ denote the cotangent line at $0$.
Writing a point $p\in P^k\setminus A^{k-1}$ in the form~(10)
we associate to this point the principal part $u/x$ with pole of
order~1. It is determined uniquely up to the action~(11)
of $\Z/k\Z$ on $u$.
This principal part determines the tangent vector at~$0$
as the linear functional on $L$:
$$
\omega\mapsto \Res_{x=0}\ \frac{u}x\omega,
$$
$\omega$ being the germ of a $1$-form representing
a cotangent vector. The $k$th tensor power of the constructed
tangent vector is an element in $(L^\vee)^{\otimes k}$,
which depends neither on the choice of the coordinate~$x$,
nor on the choice of $u$. We denote this element
by $\phi(p)\in (L^\vee)^{\otimes k}$. The mapping $\phi$
is extended continuously
to the entire $P^k$ by setting it identically~$0$
on $A^{k-1}\subset P^k$.

\proclaim{Lemma 3.2}
The dual
mapping $\phi^*$
is a morphism of graded algebras. The zero locus of $\phi$
is $A^{k-1}$. The multiplicity of $\phi$ along
$A^{k-1}$ equals $k$.
\endproclaim

\demo{Proof} For a coordinate $x$ fixed, the mapping $\phi$
is written as the polynomial $u^k$.
By definition of $P^k$, it is a well-defined polynomial
of weight~$1$ on~$P^k$.
This polynomial vanishes precisely on $A^{k-1}$.
In order to compute the multiplicity of $\phi$ along
$A^{k-1}$
let us count the number of preimages of $\phi$
that glue together at a point $a=(a_1,\dots,a_{k-1})\in A^{k-1}$
as their image $u^k$ tends to zero. These are
the principal parts with the coordinates
$(u,a_1,\dots,a_{k-1})$, $(u,\zeta a_1,\dots,\zeta^{k-1} a_{k-1})$,
$(u,\zeta^2 a_1,\dots,\zeta^{2(k-1)}a_{k-1})$, \dots,
$(u,\zeta^{k-1} a_1,\dots,\zeta^{(k-1)(k-1)}a_{k-1})$.
For a generic $a\in A^{k-1}$ all these $k$ points are distinct,
and the proof is completed.
\qed
\enddemo

\subheading{3.4. The cones of generalized principal parts}
Denote by $\Pp_i$ the sheaf of graded algebras over $\M_{g;n}$
whose stalk at a point  $(C;x_1,\dots,x_n)\in\M_{g;n}$
is the algebra of $\Z/{k_i}\Z$-invariant polynomials in
twisted Laurent coefficients of order $k_i$ at $x_i$.
Let $\PP_i=\Spec(\Pp_i)$ be the corresponding cone over
$\M_{g;n}$. The sheaf $\Aa_i$ and the corresponding cone
$\A_i=\Spec(\Aa_i)$ of generalized principal parts
with zero leading coefficients are introduced similarly.
Evidently, $\A_i$ is a subcone in $\P_i$.

The direct sum $\PP=\PP_1\oplus\dots\oplus\PP_n$ of
cones will be referred to as the
{\it cone of generalized principal parts\/}
over $\M_{g;n}$.
The cone $\PP$ contains the subcone $\A=\A_1\oplus\dots\oplus\A_n$;
the defining sheaves are denoted by $\Pp$ and
$\Aa$ respectively.
Similarly, we consider the vector bundles
$\widetilde \A_i$ and $\widetilde \A=\widetilde \A_1\oplus\dots
\oplus\widetilde \A_n$ over $\M_{g;n}$.

With each cone $\P_i$ we associate the cone
morphism $\varphi_i\:\P_i\to(\L_i^\vee)^{\otimes k_i}$,
which takes a principal part of order $k_i$
at the $i$th marked point to the $k_i$th power
of a tangent vector at this point.
The direct sum of these morphisms determines
the morphism $\varphi\:\P\to\L$, where
$\L=(\L_1^\vee)^{\otimes k_1}\oplus\dots\oplus(\L_n^\vee)^{\otimes k_n}$.

\proclaim{Lemma 3.3}
{\rm 1)}
The cone $\A_i$ is the quotient of the
vector bundle $\widetilde\A_i$ modulo the fiberwise action of the
group $\Z/k_i\Z$. The cone $\A$ is the quotient of the
vector bundle $\widetilde\A$
modulo the fiberwise action of the
group $\Z/k_1\Z\oplus\dots\oplus\Z/k_n\Z$.

{\rm 2)}
The sheaves $\Aa_i,\Aa$ are constant sheaves over $\M_{g;n}$.

{\rm 3)} The zero locus of the
morphism $\varphi_i\:\P_i\to\L_i^{\vee\otimes k_i}$
is the cone $\A_i\subset\P_i$. The multiplicity
of $\varphi_i$ along $\A_i$ equals~$k_i$.
The zero locus of the
morphism $\varphi\:\P\to\L$ is the cone $\A\subset\P$.
The multiplicity of $\varphi$ along $\A$ equals~$k_1\dots k_n$.
\endproclaim

\demo{Proof}
We need to prove only that $\Aa_i$ are constant sheaves.
This follows from assertion~4) of Lemma~3.1:
treating $a_j$, for $j=1,\dots,k_i-1$,
as coordinates along fibers of $\A_i$ and
setting $a_j=1$, $a_l=0$ for $l\ne j$ we obtain
$k_i-1$ sections of the cone $\A_i$
independent at each point of the base.
Assertion~3) of the lemma follows immediately from Lemma~3.2.
\qed
\enddemo

\subheading{3.5. Hodge bundle}
The {\it Hodge bundle\/} $\W_{g;n}\to\M_{g;n}$ is the rank $g$ vector bundle
over $\M_{g;n}$ whose fiber at a point $(C;x_1,\dots,x_n)$
is the space of holomorphic sections of the dualizing sheaf over $C$.
For a smooth $C$, this space simply coincides with the space
of holomorphic 1-forms over $C$. For a
singular curve~$C$, an element of the fiber of $\W_{g;n}$ is
a meromorphic 1-form $\omega$ on $C$ admitting poles of order at most~1
at the double points and no other poles and such that
the sum of its residues along the two branches of the curve
meeting at each double point is zero.
Below, we use the term a {\it generalized holomorphic $1$-form\/} as a synonym
for an {\it element of the Hodge bundle}.

Let $(p_1,\dots,p_n)$ be an $n$-tuple of principal parts
with poles of order $k_i$ at the marked points $x_i$
of a stable curve $(C;x_1,\dots,x_n)\in\M_{g;n}$.
Each principal part $p_i$ determines the linear functional
$$
\omega\mapsto\Res_{x_i} p_i \omega
$$
on the space of generalized holomorphic 1-forms,
given by the residue. The residue mapping naturally extends
to the entire fiber of $\P_i$ by setting the residue
equal to $0$ on the fiber of $\A_i$.

The following statement is classical, see e.g.~\cite{HMo}.

\proclaim{Theorem 3.4}
An $n$-tuple $(p_1,\dots,p_n)$ is the $n$-tuple of
generalized principal parts of
a meromorphic function $f$ over $C$ with poles of order $k_i$
at $x_i$ and no other poles if and only if
$$
\Res_{x_1} p_1\omega+\dots+\Res_{x_n} p_n\omega=0 \tag13
$$
for any generalized holomorphic 1-form $\omega$ on $C$.
Two meromorphic functions
with coinciding $n$-tuples of principal parts at the marked points
differ by an additive constant.
\endproclaim

\subheading{3.6. Completed Hurwitz spaces as cones over moduli 
spaces of curves}
Consider the mapping $R\:\PP\to\W_{g;n}^\vee$
taking an $n$-tuple $(p_1,\dots,p_n)$ of principal parts to
the linear functional $\omega\mapsto
\Res_{x_1} p_1\omega+\dots+\Res_{x_n} p_n\omega$ on $\W_{g;n}$.
Let $\ZZ$ be the zero subscheme of~$R$ in $\PP$,
i.e., $\ZZ$ is given by equations~(13).

Another way to understand $\ZZ$ is as follows.
Consider a meromorphic function on a
stable curve $C$ with poles of order $k_i$ at the marked points $x_i$
and no other poles. Associating with $f$ the $n$-tuple
of its generalized principal parts
at the marked points $x_1,\dots,x_n$ we obtain a point in $\PP$.
Now $\ZZ$ is the closure of the image of this mapping.

\remark{Definition 3.5}
The {\it completed Hurwitz space\/}
$\H_{g;k_1,\dots,k_n}$
is the closure in $\PP$ of the space of $n$-tuples of
principal parts
corresponding to meromorphic functions on smooth curves.
\endremark

The set of points of $\H_{g;k_1,\dots,k_n}$ that are not obtained from
meromorphic functions on smooth curves will be called
the {\it boundary\/} of $\H_{g;k_1,\dots,k_n}$, and denoted
by $\partial \H_{g;k_1,\dots,k_n}$.

\proclaim{Lemma 3.6}
The completed Hurwitz space has a natural structure
of a cone over $\M_{g;n}$.
\endproclaim

\demo{Proof} Consider the sheaf of ideals defining
$\H_{g;k_g,\dots,k_n}$ as a subscheme of $\PP$. Over each smooth curve
such an ideal is generated by the left hand sides of 
equations~(13),
which are homogeneous of weight~$1$, and hence the ideal is
$\C^*$-invariant. Therefore, it is  $\C^*$-invariant over each
stable curve as well. Now, the sheaf of graded algebras defining
 $\H_{g;k_g,\dots,k_n}$ as a cone is just $\PP$ modulo this sheaf of
 ideals.
\qed
\enddemo

\remark{Remark 3.7}
The completed Hurwitz space can be different from the entire
subvariety $\ZZ\subset\PP$ given by equations~(13).
In other words, not each meromorphic function on a singular
curve is the limit of a family of meromorphic functions
on smooth curves.

The simplest
example arises for the case of two poles of order~1 on elliptic
curves. Smooth elliptic curves degenerate into singular curves
having a smooth elliptic component and a smooth rational component
intersecting each other at a double point.
The stability condition implies that both marked points
belong to the rational component. For a smooth curve $C$,
the space of pairs $(p_1,p_2)$ of principal parts
of meromorphic functions at the marked points is 1-dimensional
because there is a restriction on the residues given
by the holomorphic 1-form on $C$. For a singular curve $C$,
there is no such restriction since each holomorphic 1-form
is trivial on the rational component, and the space of principal
parts of meromorphic functions is 2-dimensional. However,
not each such function can be obtained as the limit of a
family of meromorphic functions
on smooth curves. The limit functions form a 1-dimensional
subspace in the space of all meromorphic functions
with given poles of order one on the
rational component.
\endremark

\heading {4. The Lyashko--Looijenga mapping and its extension}
\endheading

\subheading{4.1. Functions with vanishing poles}
In this section we describe a possible way of thinking
about functions with vanishing principal parts.
Let $(p_1,\dots,p_n)\in \ZZ$ be a tuple of generalized principal parts,
and let $(C;x_1,\dots,x_n)\in\M_{g;n}$ be the underlying stable
curve. If $(p_1,\dots, p_n)$ are nonvanishing, then
there exists a meromorphic function $f$
on $C$ with given principal parts $p_i$ at the marked points.
This function is determined uniquely up to an additive constant.
However, if there is a generalized principal part $p_i$
with a vanishing leading coefficient, $p_i\in\A_i$,
then, for $k_i>1$, it cannot be reconstructed
from the function $f$, since the principal part of $f$
at $x_i$ is~$0$. The following construction allows one
to avoid this ambiguity.

We associate with the curve $C$ and the given tuple of generalized principal
parts another curve
$\widetilde C$ with $n$ marked points and
a meromorphic function $\tilde f$
on $\widetilde C$ determined uniquely
up to an additive constant.
If neither of the generalized principal parts $p_i$
has zero leading coefficient, then $\widetilde C$ simply coincides
with $C$, and $\tilde f$ is the same as $f$.
In general case, $\widetilde C$ is obtained from $C$ by adding
rational (i.e., of genus~$0$) components, one for each generalized
principal part with zero leading coefficient. The $i$th additional rational
component
intersects $C$ transversally at the $i$th marked point.
The new double point becomes unmarked, and a new
marked point $x_i$ distinct from the double point
is added on the new rational component. To do this we introduce a new
local coordinate $y=u/x$ and set $x_i$ to be the 
point $y=\infty$.

Note that the curve $\widetilde C$ thus constructed is no more
stable since each new rational component contains only two
special points: the double point and the marked point. However, recall
that we are not adding this unstable curve to $\M_{g;n}$, but rather
provide a better understanding of generalized principal parts with
zero leading coefficients over stable curves.

The function $\tilde f$ coincides with $f$ on $C\subset\widetilde C$,
and it is extended to the new rational components
in the following way.
Let us define the monic polynomial
$$
a(y)=y^k+a_1y^{k-1}+\dots+a_{k-1}y.
$$
According to assertion~4) of Lemma~3.1,
the coefficients of the polynomial $a(y)$
are well-defined up to the
action~(11) of the group $\Z/k\Z$, or,
what is the same, up to the change
$y\mapsto\zeta y$ of the coordinate
$y$.  The extension of $\tilde f$
to the new rational component coincides, up to
an additive constant, with this polynomial.
The constant is chosen so that to make $\tilde f$
continuous, i.e., it coincides with the value $f(x_i)$ on $C$.

The following statement is now obvious.

\proclaim{Lemma 4.1}
There exists a one-to one correspondence between
the set $\ZZ$ and the set of meromorphic functions
on curves $\widetilde C$ with poles of order $k_i$
at the $i$th marked point taken up to an additive constant.
\endproclaim

\subheading{4.2. Extension of the LL mapping}
Recall that the Lyashko--Looijenga mapping
is defined on a dense subset  $H_{g;k_1,\dots,k_n}$
in $\H_{g;k_1,\dots,k_n}$, and takes it
to the space $\C^{\mu-1}$  of monic polynomials in one variable of the form
$$
t^\mu+d_2t^{\mu-2}+\dots+d_\mu. \tag14
$$
In this section we explain how to extend the LL mapping continuously
to the completed Hurwitz space $\H_{g;k_1,\dots,k_n}$, and, moreover,
to the whole $\ZZ$.

\proclaim{Theorem 4.2}
The Lyashko--Looijenga mapping
is a finite covering outside the space of functions
with coinciding critical values and it
extends to a continuous mapping
$$
\LL\:\H_{g;k_1,\dots,k_n}\to\C^{\mu-1}
$$
equivariant with respect to the action of $\C^*$
on the left and on the right.
\endproclaim

\demo{Proof}
The assertion that the LL mapping is a covering
over the subset of polynomials with distinct roots
is a well-known corollary of the Riemann theorem
(see, e.g.,~\cite{A1}).

Let $(p_1,\dots,p_n)\in\partial\H_{g;k_1,\dots,k_n}$
be an $n$-tuple
of generalized principal parts at the marked points of a
curve $(C;x_1,\dots,x_n)$, which may be singular
and have more than one irreducible component.

Let $\tilde f\:\widetilde C\to\CP^1$
be the function corresponding to $(p_1,\dots,p_n)$
with the constant chosen arbitrarily.
It can happen that $\tilde f$ is constant on some
irreducible components of $C\subset\widetilde C$.
Denote by $C'$ the union of all irreducible
components of $C$, where $\tilde f$ is constant.

We {\it define\/} the set of critical values
of $\tilde f$ as the union of the following sets:
\roster
\item for each irreducible component of $\widetilde C$,
where $\tilde f$ is not a constant, the critical
values of the restriction of $\tilde f$ to this
component;

\item for each double point, the value of $\tilde f$ at this point
  taken with multiplicity $2$;

\item for each irreducible component of $\widetilde C$
of genus~$g'$
where $\tilde f$ is a constant, the value of $\tilde f$
on this component taken with multiplicity $2g'-2$.
\endroster

Now let us verify that the LL mapping thus extended
is continuous.

Take a double point $\theta$ of the curve $\widetilde C$ and
consider the preimage $\tilde f^{-1}(S^1)$ of a small circle $S^1$ in $\C$
centered at $\tilde f(\theta)$. We are interested in the intersection
of  this  preimage with a small neighborhood of $\theta$ in
$\widetilde C$. If $\tilde f$ is nonconstant
on both branches meeting at $\theta$, 
then this intersection consists of
two circles, one on each branch. If $\tilde f$ is nonconstant
only on one such branch, then there is only one circle, and it lies on
this branch. Finally, if $\tilde f$ is constant
on both branches, 
then the intersection in question is empty, and the
corresponding double point is called an {\it inner\/} double point.

Each of the circles carries a natural number, the degree
of its mapping to $S^1$ under $\tilde f$. This degree coincides
with the order of the critical point of $\tilde f$
on the corresponding branch increased by~$1$.

Consider now a holomorphic deformation $F\:\uC\to\CP^1$ of $\tilde f$,
where $\uC\to\C^1$ is a holomorphic family of
marked genus $g$ curves
such that $C_\tau$ is smooth for $\tau\ne0$
and $C_0=\widetilde C$, and the restriction $f_\tau$
of $F$ to $C_\tau$ is a meromorphic function, $f_0=\tilde f$.
We do not require that $f_\tau$ has zero sum of
critical values.

For a fixed double point $\theta$,  consider the preimage $F^{-1}(S^1)$ of
the same circle  $S^1$ in $\C$ centered at $\tilde f(\theta)$. Take
the connected component of this preimage containing one of the circles
on $\widetilde C$ built above. For
a sufficiently small $\tau\ne0$, the intersection of this connected
component with $C_\tau$ is a circle, and the degree of its mapping to
$S^1$ under $f_\tau$ coincides with the number assigned to the circle
on $\widetilde C$.

Applying the procedure described above
for all double points of $\widetilde C$
we obtain a finite set of circles on each curve $C_\tau$
for $\tau$ small enough. These circles cut each
curve $C_\tau$ into connected pieces of three
different types:
\roster
\item pieces, containing marked points. These pieces
are in one-to-one correspondence with the irreducible
components of $\widetilde C$, where $\tilde f$ is nonconstant;
\item pieces, holomorphicaly equivalent to the annulus.
These pieces are in one-to-one correspondence with those
double points of $\widetilde C$, where $\tilde f$ is nonconstant
on both branches meeting at this point;
\item  pieces of positive genus without marked points.
These pieces are in one-to-one correspondence
with the connected components of the curve $C'\subset C$,
the constant locus of $\tilde f$.
\endroster

Now let us follow the behavior of the critical points
and critical values of $f_\tau$ on pieces of each
three types. Let the index $j$ run over the pieces.

(1) By the Riemann--Hurwitz formula, the number of critical points
  of $f_\tau$ on a piece
of the first type is equal to $K_j+n_j+2g_j-2-D_j+c_j$,
where $K_j$ is the total
order of all poles on this piece, $n_j$ is the number
of marked points, $g_j$ is the genus of the corresponding component of
$\widetilde C$, $D_j$ is the
total degree of the circles bounding the piece, and $c_j$ is the number
of these circles.
As $\tau$ tends to~$0$, these critical points tend to
the critical points of $\tilde f$ on the corresponding component
of $\widetilde C$ with the double point excluded, and
the critical values tend to that
of $\tilde f$.

(2) The mapping $f_\tau$ takes an annulus
without marked points to the small disk
bounded by $S^1$. The degree of the mapping
is the sum of integers assigned
to the boundary circles. By the Riemann--Hurwitz formula,
the number of critical points
of $f_\tau$ on the annulus coincides with the degree.
As $\tau$ tends to~$0$, these critical points tend to
the double point, and the critical values tend to
the value of $\tilde f$ at the double point.

(3) The mapping $f_\tau$ takes a piece of the third
type to a disk. The degree of this mapping equals $D_j$, the
sum of numbers assigned to the boundary circles.
The number of
critical values on this piece is $2g_j-2+D_j+c_j$,
where $c_j$ is the number of boundary circles, and $g_j$ is the
(arithmetic) genus of the corresponding connected component
$C'\subset\widetilde C$.
Further, $2g_j-2=2\sum(g_{ji} -2)+2l$, where $l$ is the number of inner
double points on $C'$, and $g_{ji}$ are the genera of irreducible
components of $C'$.
As $\tau$ tends to~$0$, the critical values tend to
the constant value of $\tilde f$ on the limit curve.

Now, taking the union of all limit critical
values over all pieces we conclude
that the set of critical values of $f_\tau$
tends precisely to the set described above.
\qed
\enddemo

It is now easy to get the following Corollary.

\proclaim{Corollary 4.3}
The LL mapping extended on the whole $\ZZ$ as above is continuous.
\endproclaim

\heading {5. Top Segre classes}
\endheading

\subheading{5.1. Segre classes of vector bundles and cones}
This section follows (with slight modification)
the approach of~\cite{Fu, Ch.~3, 4, 8, 19}.

Let $X$ be  a nonsingular variety or an orbifold of pure dimension $d$.
%of rank $r+1$
Consider a pure dimensional cone $\cS$ over
$X$, and let $\pi\:P\cS\to X$
denote the projection of the corresponding weighted projective cone.

The {\it Segre class\/} of $\cS$
is an element in
the cohomology ring $H^*X=H^*(X,\Q)$ of $X$ with rational coefficients;
it is defined as follows.
Consider the value $c_1(\O(1))^{i}\cap[P\cS]$
of the iterated first Chern class
of the canonical line bundle on the fundamental cycle
of $P\cS$; here Chern classes of bundles over cones are considered in
the operational
sense (see \cite{Fu, Ch.~3}). This value is an element in the group
$A_*P\cS$ of cycles modulo the rational
equivalence on the projective
cone $P\cS$.
The mapping $\pi_*\:A_*P\cS\to A_*X$ pushes this
element to a cycle in $X$. The isomorphism
$A_*X\cong A^{d-*}X$ takes this cycle to the Chow ring
$A^*X$ (\cite{Fu, Sec.~8.3}).
In its own turn, the Chow ring admits a homomorphic
graded mapping ${\cl}\:A^*X\to H^*X$, the class mapping
(\cite{Fu, Corollary~19.2}).
Denote the composition of all these mappings
by $h\:A_* P\cS\to H^*X$. The Segre class of $\cS$ is defined
by the expression
$$
s(\cS)=h\bigg(\sum_{i\gs0}c_1(\O(1))^{i}\cap[P\cS]\bigg).
$$
It can be represented uniquely in the form
$$
s(\cS)=s_0(\cS)+\dots+s_d(\cS),\qquad s_i(\cS)\in H^{2i}X.
$$

The class $s_d(\cS)$ is called the {\it top Segre class\/}
of the cone $\cS$ and denoted by $s_{top}(\cS)$; its {\it value\/}
$\int_Xs_{top}(\cS)$
on $X$ is a rational number.
Lemma~5.3 below expresses the multiplicity of the LL mapping
in terms of the top Segre class of the cone $\H_{g;k_1,\dots,k_n}$.

For a vector bundle $E$ over $X$, the Segre class
coincides with the inverse Chern class,
$s(E)=c^{-1}(E)$; here the Chern class is understood in the usual
sense, as an element in the cohomology ring
$H^*X$.

\subheading{5.2. Segre classes of cones of principal parts}
Now we are going to compute the  Segre classes
we shall require below.

\proclaim{Lemma 5.1}
The Segre class of the cone $\A$ equals
$$
s(\A)=\prod_{i=1}^n\frac1{k_i}\frac{k_i^{k_i-1}}{(k_i-1)!}.
$$
\endproclaim

\demo{Proof}
Since, by Lemma~3.3, the sheaf $\Aa_i$,
which determines the cone $\A_i$,
is a constant sheaf, the only nontrivial Segre class
of $\A_i$ is its zero Segre class $s_0(\A_i)$.
The last one can be found locally,
at a fiber of $P\A_i$ over a point in $\M_{g;n}$.
This fiber is the projectivized weighted vector space
$\widetilde A^{k_i-1}$ modulo the action of the group $\Z/k_i\Z$.
Its zero Segre class is the inverse product of
the weights of the coordinates divided by the order
of the acting group,
$$
s_0(\A_i)=\frac1{k_i}\frac{k_i^{k_i-1}}{(k_i-1)!}.
$$

 In a similar way we get that $s_0(\A)$ is the product of Segre
 classes $s_0(\A_i)$.
\qed
\enddemo

\proclaim{Lemma 5.2}
The Segre class of the cone $\PP$ equals
$$
s(\PP)=\prod_{i=1}^n\frac{k_i^{k_i}}{k_i!}\frac{1}{1-k_i\psi_i},
$$
where $\psi_i=c_1(\L_i)$.
\endproclaim

\demo{Proof}
Consider the  sequence
of cone morphisms
$$
0\longrightarrow
\A\longrightarrow\PP\overset\varphi\to\longrightarrow\L
\longrightarrow0,
$$
where $\A\to\PP$ is the embedding, and $\varphi\:\PP\to\L$
is as in Section~3.4.

Let $\pi\:P\PP\to\M_{g;n}$ be the projection of
the projective cone of principal parts to the
base. This projection defines the bundle $\O(1)\otimes\pi^*(\L)\to P\PP$.
The mapping $\varphi$ determines a section
$\sigma_\varphi\:P\PP\to \O(1)\otimes\pi^*(\L)$ of
this bundle in the following way. 
The value of $\sigma_\varphi$ at a point
$\C^*p\in P\P$ is equal to $1_p\otimes\pi^*(\varphi(p))$;
here $1_p$ is the element in the fiber of~$\O(1)$
over $\C^*p$ that equals~$1$ on~$p$.

The zero locus $Z(\sigma_\varphi)\subset P\PP$
of the section $\sigma_\varphi$ is the
projectivization of the zero locus of~$\varphi$.
By Lemma~3.3, it is $P\A\subset P\PP$
taken with multiplicity $k_1\dots k_n$.
It is of pure codimension $n$ in $P\P$. Therefore,
by~\cite{Fu, Proposition~14.1},
this zero locus represents the top Chern class
of the bundle $\O(1)\otimes\pi^*(\L)\to P\PP$,
$$\align
[Z(\sigma_\varphi)]&=c_n(O(1)\otimes\pi^*(\L)) \cap[P\PP]\\
&=\sum_j c_{n-j}(\pi^*(\L))c_1(\O(1))^j\cap [P\PP].
\endalign
$$
Thus,
$$\align
k_1\dots k_ns(\A)&=\sum_lc_1(\O(1))^l\cap[Z(\sigma_\varphi)]\\
&=\sum_{l,j} c_{n-j}(\pi^*(\L))c_1(\O(1))^{l+j}\cap[P\PP]
\\
&=c(\L)s(\PP).
\endalign
$$
Hence,
$$
s(\PP)=k_1\dots k_n\frac{s(\A)}{c(\L)},
$$
and the required assertion follows.
\qed
\enddemo

\subheading{5.3. Multiplicity of the LL mapping and Segre classes}
Similarly to the proof of Lemma~5.2, one can associate
to the
mapping $R\:\PP\to\W_{g;n}^\vee$ given by the sum of residues
at the marked points the section
$$
\sigma_R\:P\PP\to\O(1)\otimes\pi^*(\W_{g;n}^\vee)
$$
of the vector bundle over the projective cone $P\PP$:
$\sigma_R$ associates to a point in $P\PP$
the linear function on the fiber of $\pi^*(\W_{g;n})$
determined by the residue with the element in
the fiber of $\O(1)$. By Theorem~3.4,
the zero locus of the section $\sigma_R$
is the projectivization $P\ZZ$ of the space
$\ZZ$ of meromorphic functions.

Recall that the completed Hurwitz space $\H_{g;k_1,\dots,k_n}$
is a subcone in the cone $\PP$ of principal parts.

\proclaim{Lemma 5.3}
The multiplicity of the LL mapping
$\LL:\H_{g;k_1,\dots,k_n}\to\C^{\mu-1}$
is equal to
$$
\mu!\int_{\M_{g;n}}s(\H_{g;k_1,\dots,k_n}).
$$
\endproclaim

\demo{Proof}
Denote by $\Discr\subset\C^{\mu-1}$ the discriminant in
the target space of the LL mapping, i.e., the subvariety
of polynomials having at least two coinciding roots.
On the Hurwitz space, the LL mapping
is a local covering, and the image of
the boundary $\partial\H_{g;k_1,\dots,k_n}$
is of codimension at least one in the target space.
Indeed, for any element of  $\partial\H_{g;k_1,\dots,k_n}$
the underlying curve is singular,
and on a singular curve at least two critical
values of a meromorphic function coincide, i.e.,
$\LL(\partial\H_{g;k_1,\dots,k_n})\subset\Discr$.
Hence, in order to compute the multiplicity of LL
it suffices to compute the number of preimages
of a generic point in $\C^{\mu-1}$ or, what is the same,
the number of $\C^*$-orbits in the preimage of a $\C^*$-orbit.

Define a {\it standard hypersurface\/} in the space~$\C^{\mu-1}$
of polynomials as the image of a hyperplane
in the space of roots with zero sum under the Vieta mapping
$V\:\C^{\mu-1}\to\C^{\mu-1}$ taking the set of roots to the set of
coefficients.
A generic orbit of the $\C^*$-action in the space
of coefficients is an irreducible component in the intersection of~$\mu-2$
standard hypersurfaces in general position.
Therefore, in order to compute the multiplicity of the LL mapping
it is sufficient to compute the selfintersection index
of the preimage of a standard hypersurface
under LL.

The preimage of a standard hypersurface under the Vieta mapping
consists  of $\mu!$
hyperplanes, hence the preimage of a standard hypersurface
under the LL mapping corresponds to $\mu!c_1(\O(1))$.
The intersection index of
$\mu-2$ such preimages is precisely the degree of
$$
(\mu!)^{\mu-2}(c_1(\O(1))^{\mu-2}\cap [P\H_{g;k_1,\dots,k_n}]),
$$
or, which is the same, the value of
$s_{top}(\H_{g;k_1,\dots,k_n})$ on
$\M_{g;n}$ times $(\mu!)^{\mu-2}$.

Now consider the intersection of preimages
of $\mu-2$ standard hypersurfaces under the LL mapping. Evidently, the
number of irreducible components in
this intersection divided by the multiplicity of the LL mapping
equals the number of intersections of preimages of the same
standard hypersurfaces under the Vieta mapping divided by the multiplicity
of the Vieta mapping. However, in the space of roots the preimages of
standard hypersurfaces are just the sets of hyperplanes, so 
in this space the intersection contains $(\mu!)^{\mu-2}$
irreducible components. Recalling that the multiplicity
of the Vieta mapping is $\mu!$ we obtain the required expression.
\qed
\enddemo

By the definition of the LL mapping,
$\LL^{-1}(\C^{\mu-1}\setminus\Discr)\subset\H_{g;k_1,\dots,k_n}$.
Denote by $\ZZ'\subset\ZZ$ the union of all irreducible
components of $\ZZ$ that are not contained in $\H_{g;k_1,\dots,k_n}$.
We have $\ZZ=\H_{g;k_1,\dots,k_n}\cup\ZZ'$ and
$\LL(\ZZ')\subset\Discr$.

\proclaim{Lemma 5.4}
$c_1(\O(1))^{\mu-2}\cap[P\ZZ']=0$.
\endproclaim

\demo{Proof}
The inverse image $\LL^{-1}(0)$ is the zero section of $\ZZ$
(all twisted Laurent coefficients equal $0$). Therefore,
the projectivized LL mapping $P\LL$ 
is well-defined and it maps $P\ZZ$ to
$P\C^{\mu-1}$ and $\O_{P\ZZ}(1)=P\LL^*\O_{P\C^{\mu-1}}(1)$.
In particular, $\O_{P\ZZ'}(1)=(P\LL|_{P\ZZ'})^*\O_{P\Discr}(1)$.
But, since
$\dim\ P\Discr=\mu-3$,
the bundle $\O(1)$ over $P\ZZ'$ is induced from a  line
bundle on a variety of dimension
$\mu-3$, and the lemma is proved.
\qed
\enddemo

\proclaim{Lemma 5.5}
$$
\int_{\M_{g;n}}s(\H_{g;k_1,\dots,k_n})=
\int_{\M_{g;n}}c(\W_{g;n}^\vee)s(\PP).
$$
\endproclaim

\demo{Proof} Recall that $P\ZZ\subset P\PP$ is the zero locus
of the section $\sigma_R\:P\PP\to\O(1)\otimes\pi^*(\W_{g;n}^\vee)$
of the bundle $\O(1)\otimes\pi^*(\W_{g;n}^\vee)$ over $ P\PP$.
We start with constructing a $(\mu-2)$-cycle $D\subset P\ZZ$,
$[D]\in A_{\mu-2}P\ZZ$ representing the localized top Chern class
of the bundle $\O(1)\otimes\pi^*(\W_{g;n}^\vee)$ over $ P\PP$,
$[D]=c_g (\O(1)\otimes\pi^*(\W_{g;n}^\vee))\cap [P\PP]$.

Such a construction works for an arbitrary vector bundle $E\to X$
over a pure dimensional scheme $X$ 
and a section $\eta\:X\to E$ of this bundle
(see~\cite{Fu, Secs.~14.1, 6.1}).
Let $e$ denote the rank of $E$, and let $Z(\eta)\subset X$
be the zero locus of $\eta$. The normal cone $N$
to $Z(\eta)$ in $X$
is naturally embedded, as a subcone, in the total space of the
vector bundle $E$ restricted to $Z(\eta)$ in the following way.
A tangent vector $\xi\in T_z X$ to $X$ at a point $z\in Z(\eta)$
is taken to the tangent vector
$d\eta(\xi)\in T_{\eta(z)}E$ treated as an element
of the fiber $E_z$ of $E$ at $z$, which is identified
naturally with the quotient space $E_z=T_{\eta(z)} E/T_z X$.
The normal cone to each irreducible component of $Z(\eta)$ in
$X$ has the same dimension as $X$ itself,
and $[N]$ is a
cycle of codimension $e$ in $A_*E|_{Z(\eta)}$.
Proposition~14.1 from~\cite{Fu} states that the image of $N$
under the isomorphism $\tau\: A_*E\cong A_{*-e}X$ is
a codimension $e$ cycle in $A_{*-e}X$ representing the localized
top Chern class of the bundle $E$ on $X$. 

In particular, each irreducible component of $Z(\eta)$
of codimension~$e$ in~$X$ enters the class~$[D]$
with some positive multiplicity, which is determined
by the behavior of the section~$\eta$ along this component.

Applying the construction above to the vector bundle
$\O(1)\otimes\pi^*(\W_{g;n}^\vee)$
over $P\PP$ and the section $\sigma_R$ we obtain a
cycle $[D]\in A_{\mu-2}P\PP$ such that
$$
[D]=c_g(\O(1)\otimes\pi^*(\W_{g;n}^\vee))\cap [P\PP]=
\sum_ic_1(\O(1))^ic_{g-i}(\pi^*(\W_{g;n}^\vee))\cap [P\PP].
$$
Therefore, one has
$$\align
\pi_*\bigg(\sum_j c_1(\O(1))^{j}&\cap[D]\bigg)\\
&=\pi_*\bigg(\sum_jc_1(\O(1))^{j}\sum_ic_1(\O(1))^i
c_{g-i}(\pi^*(\W_{g;n}^\vee))\cap[P\PP]\bigg)\\
&=\pi_*\bigg(\sum_{i,j}c_1(\O(1))^{i+j}c_{g-i}
(\pi^*(\W_{g;n}^\vee))\cap[P\PP]\bigg)\\
&=\pi_*\bigg(\sum_i c_i(\pi^*(\W_{g;n}^\vee))\sum_{l=i+j}
c_1(\O(1))^{l}\cap[P\ZZ]\bigg).
\endalign
$$
Hence $s(D)=c(\W_{g;n}^\vee)s(\PP)$; in particular,
$$
\int_{\M_{g;n}}s(D)=\int_{\M_{g;n}}c(\W_{g;n}^\vee)s(\PP).
$$

The construction above provides a representative $D$ 
of the class $[D]$ that is contained in $P\ZZ$. 
We set $[D_\H]=[D\cap P\H_{g;k_1,\dots,k_n}]$
and $[D']=[D\cap P\ZZ']$ so that $[D]=[D_\H]\cup[D']$.
%The cycle $[D]$ can be written as the sum $[D]=[D_\H]\cup[D']$,
%where $[D_\H]=[D]\cap P\H_{g;k_1,\dots,k_n}$ and $[D']=[D]\cap\ZZ'$.
%By Lemma~5.4,
%$c_1(\O(1))^{\mu-2}\cap [D']=0$, and 
%the only thing we need to prove is that 
Let us prove now that $P\H_{g;k_1,\dots,k_n}$ 
represents the class $[D_\H]$. As we have seen above, this 
is true up to multiplicity, since the codimension 
of $P\H_{g;k_1,\dots,k_n}$ in $P\PP$ is~$g$. 
It remains to prove that the section 
$\sigma_R$   
is transversal to the zero section of the bundle 
$\O(1)\otimes\pi^*(\W_{g;n}^\vee)$
at a generic point of $P\H_{g;k_1,\dots,k_n}$,
and hence the multiplicity in question is precisely~$1$.

In the Appendix we give a proof of this statement in a more general 
situation. Here we try to explain it without using sophisticated tools
from algebraic geometry.

Fix a generic point in $P\H_{g;k_1,\dots,k_n}$.
Such a point is a smooth stable curve
$(C^0;x^0_1,\dots,x^0_n)$ endowed with an $n$-tuple
$(p^0_1,\dots,p^0_n)$ of generalized principal parts
at the marked points taken up to a common
nonzero multiplier and such that:

1) the leading terms of the principal parts do not vanish;

2) the mapping $f^0\:C^0\to\CP^1$ with given principal
parts has distinct critical values;

3) the principal parts $p^0_i$ satisfy equations~(13)
for any holomorphic 1-form $\omega$ on $C^0$.

For the sake of brevity, below we denote the chosen point
simply by $(C^0,f^0)$.

The tangent space to the total space of the bundle
$\O(1)\otimes\pi^*(\Lambda^\vee_{g;n})$
over $P\PP$ at $(C^0,f^0)$
splits naturally into the direct sum of the horizontal
tangent space along the base and the vertical tangent
space along the fiber of the bundle.
The vertical tangent space is naturally identified
with the space $\Lambda^\vee(C^0)$ of
linear functionals on the space of holomorphic
1-forms on $C^0$.
The section $\sigma_R$ is transversal to the zero section 
at the point $(C^0,f^0)$ if and only if
the projection of the image of $d\sigma_R$ at this point
to the vertical tangent space is a surjection.
We denote the composition of $d\sigma_R$ with
the projection to the vertical tangent space by
$d\tilde\sigma_R\:T_{(C^0,f^0)}P\P\to \Lambda^\vee(C^0)$.

Below, we distinguish between two cases. In each case, we are going
to construct a subspace $L$ of the tangent space $T_{(C^0,f^0)}P\P$
such that the restriction of $d\tilde\sigma_R$ onto $L$ is of rank $g$.

1. $k_1+\dots+k_n>g$. We choose for $L$ the tangent space to the fiber of
$\pi$; that is,  
$L=T_{(C^0,f^0)}\pi^{-1}(C^0;x_1^0,\dots,x_n^0)$. Then the restriction
of $d\tilde\sigma_R$ onto $L$ is defined by the $g\times K$ Brill--Noether
matrix (see \cite{ACGH, Ch.~4}). Since in this case $C^0$ is a generic curve
and $f^0$ is a generic function, the rank of the Brill--Noether matrix
equals $g$.

2. $k_1+\dots+k_n\ls g$. This situation can happen only if
$g\gs2$ since there are no meromorphic functions on
elliptic curves with a single pole of order~1.

We start with associating to each
noncritical unmarked point $y\in C^0$
the germ of a holomorphic curve in $P\PP$ passing through
$(C^0,f^0)$. For each such point $y\in C^0$,
the function $f^0$ determines a coordinate $z_y=f^0-f^0(y)$
in a neighborhood of $y$.

Fix the complex structure of $C^0$ outside a small disc $U_{2r}\subset C^0$
of radius $2r\ls2$ (in the coordinate $z_y$)
centered at $y$ and not containing
any of the critical points of $f^0$.
Denote by $B\subset\C^1$ the standard unit disc
with the coordinate $w$ and let $W_\delta\subset B$
be the varying annulus $|w|>|\delta|^{1/2}$.
Map $W_\delta$ to $U_{2r}$ by setting $z_y=r\left(w+\frac\delta{w}\right)$.
This amounts to identifying the annulus $W_\delta$ with its image obtaining 
thus a new complex curve, which we denote by $C^\delta$.
The tangent vector to the family $C^\delta$
of deformations of the complex structure thus constructed 
is called the {\it first-order
Schiffer variation\/} centered at the point $y$,
see~\cite{HMo, Ch.~3B}. It can be constructed
starting from an arbitrary local coordinate $z_y$,
and it does not depend, up to a constant
factor, on the chosen coordinate.

Observe that the Schiffer variation centered at $y$ preserves the curve 
$C^0$ outside $U_{2r}$. Let us fix the positions of the marked points, the 
corresponding local coordinates in their neighborhoods, and the Laurent 
polynomials defining principal parts with respect to these local coordinates.
Then the first-order Schiffer variation lifts
to a tangent vector $\tau_y\in T_{(C^0,f^0)}P\P$.

Let us prove that the mapping
$d\tilde\sigma_R$
takes the vector $\tau_y$ to the following functional
on the space of holomorphic 1-forms on $C^0$ 
closely related to the canonical mapping:
$$
\omega\mapsto\Res_y\frac{\omega}{z_y}.
$$
Indeed, let  us fix an arbitrary holomorphic connection on the bundle
$\Lambda_{g;n}$ restricted to $C^\delta$, say
the one preserving the integrals
of holomorphic 1-forms over $a$-cycles.
Then for each holomorphic 1-form $\omega^0$ on $C^0$
we have a well-defined holomorphic family $\omega^\delta$
of holomorphic 1-forms over $C^\delta$. 

The 1-form $\omega^\delta$ written in the coordinate $w$ looks like
$$
\omega^\delta=(c_0+c_1w+c_2w^2+\dots)dw+\delta\eta +o(\delta),
$$
where $\eta$ is a holomorphic 1-form in the unit disc. Rewriting it in the 
coordinate $z_y$ in the image of $W_\delta$ we get
$$
\omega^\delta=\left(\frac{c_0}r+\frac{c_1}{r^2}z_y+\frac{c_2}{r^3}z_y^2+
\dots\right)dz_y+\frac{\delta c_0r}{z_y^2}dz_y+\delta\bar\eta+o(\delta),
$$
where $\bar\eta$ is a holomorphic 1-form in $U_{2r}$. Therefore,
$$
\omega^0=\left(\frac{c_0}r+\frac{c_1}{r^2}z_y+\frac{c_2}{r^3}z_y^2+
\dots\right)dz_y
$$ 
and $c_0=r\Res_y\frac{\omega^0}{z_y}$. Finally, the mapping $d\tilde\sigma_R$
takes $\tau_y$ to 
$$
\lim_{\delta\to 0}\frac1\delta\sum_i\Res_{x_i}p_i\omega^\delta=
-\lim_{\delta\to 0}\frac1\delta\int_{\partial
U_{2r}}f^0\omega^\delta=-c_0r,
$$
and the required claim follows.

Choosing points $y_1,\dots,y_g\in C^0$ in general position
we obtain $g$ vectors $\tau_{y_1},\dots,\tau_{y_g}$
whose images under the mapping
$d\tilde\sigma_R$ are linearly independent, since
the images of $g$ points on $C^0$
in general position under
the canonical mapping span a $(g-1)$-dimensional
subspace.

Hence
$$
c_1^{\mu-2}(\O(1))\cap [D]= c_1^{\mu-2}(\O(1))\cap [P\H_{g;k_1,\dots,k_n}]+
c_1^{\mu-2}(\O(1))\cap [D'].
$$

By Lemma~5.4,
$c_1(\O(1))^{\mu-2}\cap [D']=0$, and hence
$$
\pi_*(c_1(\O(1))^{\mu-2}\cap[P\H_{g;k_1,\dots,k_n}])=
\pi_*(c_1(\O(1))^{\mu-2}\cap[D]),
$$
and the result follows.
\qed
\enddemo

Combining the assertions of Lemmas~5.2, 5.3, and~5.5,
we obtain the proof of the main theorem.

\heading {Appendix: Local study of the LL map}
\endheading

Our setup will be the following. We have a smooth and proper map $f\: C\to S$ 
of relative dimension $1$. Even though many of the results will be true 
generally, we assume that $\cQ$ is contained in $\Gamma(S,\sO_S)$. 
We fix a genus $g$  and a
degree $d$ and consider the moduli stack $\sM_{g;d}(C/S)$ of families of
semistable curves $D \overset{k}\to\to T$ over an $S$-scheme $T$ of genus $g$
together with an $S$-map $D \overset{h}\to\to C$ fiberwise finite and of 
degree $d$. Our aim is first to define a stratification by subschemes of 
$\sM_{g;d}(C/S)$
on whose strata the local structure of the map $D \to S$ is fixed. For this we
will need the following lemma.

\proclaim{Definition--Lemma A.1}
Let $\pi\: X\to S$ be a finite flat map, where $S$ is a scheme over $\Sp \cQ$.

{\rm (i)}
If the map $\pi_*\sO_X \to \Hom_{\sO_S}(\pi_*\sO_X,\sO_S)$ 
induced by the relative
trace map has locally constant rank {\rm(}in the schematic sense{\rm)}, 
then there exists a {\rm (}unique{\rm)} closed subscheme 
$X^{fred} \hookrightarrow X$ such that for any
geometric point $\bar s \to S$, $X^{fred}_{\bar s}$ is the reduced subscheme of
$X_{\bar s}$. Furthermore, $X^{fred}$ is finite \'etale over $S$. We will call
this subscheme the {\rm fiberwise reduced subscheme}.

{\rm (ii)}
If $S$ is reduced, then
$\pi_*\sO_X \to\Hom_{\sO_S}(\pi_*\sO_X,\sO_S)$  is of locally constant rank 
if and only if the function
that to a point of $S$ associates the number of geometric points 
above it on $X$
is locally constant.
\endproclaim

\demo{Proof}
Consider the relative trace map $\Tr\:{\pi_*\sO_X \to \sO_S}$ and the 
associated
trace pairing $(f,g) \mapsto \Tr(fg)$. This induces a map $\pi_*\sO_X \to
\Hom_{\sO_S}(\pi_*\sO_X,\sO_S)$. By assumption, the codimension of this 
radical 
is locally constant, so that the trace pairing has locally constant rank and 
the kernel of
$\pi_*\sO_X \to \Hom_{\sO_S}(\pi_*\sO_X,\sO_S)$ is a sub-vector bundle. 
This kernel
is clearly an ideal, and the closed subscheme defined by it is 
fiberwise
reduced subscheme. As $\pi_*\sO_X \to \Hom_{\sO_S}(\pi_*\sO_X,\sO_S)$ has 
locally
constant rank, this subscheme is flat; it is also clearly fiberwise \'etale, 
and hence \'etale, and it is clearly finite.

Finally, as we are in characteristic zero, for any geometric point $\bar s \to
S$, the kernel of the map $\pi_*\sO_X \to \Hom_{\sO_S}(\pi_*\sO_X,\sO_S)$ 
equals the
radical of the structure sheaf of $X_{\bar s}$. Hence, its rank is the number 
of geometric points over $\bar s$.
\qed
\enddemo

We want to apply this lemma to the ramification and branch loci, but before we
can do that we need to define these loci in the presence of singularities 
of the fibers of $D \to S$. To do this we first note that outside of the
singular locus of $D \to S$ we have a map 
$h^*\Omega^1_{C/S} \to \Omega^1_{D/S}$,
and by tensoring by the inverse of $h^*\Omega^1_{C/S}$ we get a global section
of $\omega_{D/C}$. Now this sheaf is identified, by duality, with
$\Hom_{\sO_C}(h_*\sO_D,\sO_C)$ and we claim that under this identification the
section corresponds to the trace map $h_*\sO_D \to \sO_C$. Indeed, this can be
verified on the open dense subset where $h$ is \'etale and then, being local in
the \'etale topology, we may assume that $D$ is a disjoint union of copies of
$C$ in which case we reduce to the case of one copy by the fact that both
sides are additive over disjoint unions. This latter case is obvious.

In our case we can now go backwards; as $h$ is finite and $f$ smooth, 
$h$ is flat
and it is assumed to be finite, so we may consider its trace map as a global
section of $\omega_{D/C}$, the relative dualizing sheaf, and then 
tensoring with
$h^*\omega_{C/S}$ we get a map $h^*\omega_{C/S} \to \omega_{D/S}$. 
Now these two
sheaves are $S$-flat and the map, which commutes with base change, is fiberwise
injective, and hence the map is injective with an $S$-flat cokernel. We tensor
the map with the inverse of $\omega_{D/S}$ to define a finite and $S$-flat
subscheme $R$ of $D$ by the exact sequence
$$
0\to{\omega^{-1}_{D/C}}\to{\sO_D}\to{\sO_R}\to 0.
$$
This is by definition the {\it ramification locus\/} of $h$. We then define
the {\it branch locus}, $B$, as $\Div(h_*\sO_R)$ (cf.
\cite{Mu, Ch.~5, \S 3 }). Then the branch locus also commutes with base 
change and is finite
and flat over $S$.

\remark{Remark A.2}
This definition is a special case of \cite{FnP}, the definition
given here is somewhat more explicit. A consequence of that is that 
it is clear, by a local calculation, that it extends to the case when 
the base also is allowed to have nodes but the map then is required 
to have the same local degree
on both branches when a node maps to a node (which is the condition of
\cite{HMu}).
\endremark

A stratification by locally closed subschemes of  $\sM_{g;d}(C/S)$ will of
course induce one on any family of maps. It will however be convenient to
directly define this stratification for any family. Hence we consider a smooth
proper curve $f\:C\to S$ over a base $S$ over $\Sp\cQ$, a semi-stable curve 
$k\: D\to S$ 
and a finite $S$-map $h\: D\to C$.

We now stratify $S$ by closed subschemes by the condition that the
rank of the maps $k_*\sO_R \to \Hom_{\sO_S}(k_*\sO_R,\sO_S)$ and $f_*\sO_B \to
\Hom_{\sO_S}(f_*\sO_B,\sO_S)$ induced by the trace maps is less than or equal 
to some integers. Looking at a locally closed
stratum $T$ that is the complement in one of the closed strata of the strata
contained in it, we see that the restriction of $R$ and $B$ to $T$ fulfill the
conditions of Lemma~A.1, and hence we may define $R^{fred}$
and $B^{fred}$. Clearly the map $R \to B$ induces a map $R^{fred} \to
B^{fred}$. We associate to each point of $R^{fred}$ the length of $R$ at that
point. This is a locally constant function on $R^{fred}$, 
and as $R^{fred} \to T$ is finite \'etale, 
$T$ can be written as a disjoint union of open subschemes over
each of which this length is constant on each irreducible component of
$S^{fred}$. We may similarly refine $T$ so that the same is true for $B^{fred}$
and the length of $B$. It is clear that over such a component the local degrees
of $h$ are the same. We call these components the {\it equisingular
strata\/} of $S$. On each such stratum $T$ we have associated a
relative effective Cartier divisor $B^{fred}$ of $C \to S$ of fixed degree. 
This map $T \to \Div(C/S)$ is the LL map.

\proclaim{Proposition A.3}
Let $T$ be an equisingular stratum of $\sM_{g;d}(C/S)$. Then the LL map $T \to
\Div(C/S)$ is \'etale.
\endproclaim

\demo{Proof}
To prove that it is enough to verify the infinitesimal lifting property, so we
may assume that $S = \Sp A$ where $A$ is a local artinian ring with
algebraically closed residue field, $A'$ is another artinian local ring,
that we
have a surjective map $A \to A'$, a relative effective Cartier divisor
$B$ of 
$C/S$, \'etale over $S$, and a lifting of it to a map 
$S' := \Sp A' \to T$ over $A'$, and
we want to show that there is a unique lifting over $A$.

Rearranging our data we have a semistable curve $C$ over $S$, a relative
effective Cartier divisor $B$ of $C/S$, \'etale over $S$, a curve 
$D \to S'$ and
an equisingular $S'$-map $f\: D\to{C_{|S'}}$ whose branch locus is $B_{|S'}$. 
We
then want to show that there is a unique equisingular extension of $f$ to $S$
whose branch locus is $B$. Outside of the support of $B$ the map must be an
\'etale cover, and the existence and uniqueness is clear. Hence we may work
locally in the \'etale topology around the points of $\Supp B$. Hence we may
assume that $C=\Sp A\{t\}$, where $A\{t\}$ is the strict henselisation of 
$A[t]$
at $t=0$, and that $B$ is given by $t=0$. Then the fiberwise 
reduced subscheme of
the ramification locus of $f$ is a disjoint union of copies of $C_{|S'}$ and we
may similarly restrict ourselves to one of these components, so that we may
assume that also $D$ is strict henselian and with an $S'$-section. Having
arrived at this point we may also complete $C$ along the branch locus and
hence assume that $C=\Sp A [[t]]$ and then also $D$ will be the 
spectrum of a complete local ring $T$.

Assume first that $D$ is smooth over $S'$. We begin by showing that 
the lifting,
if it exists, is unique. For this we assume that we can choose a generator $s'$
for the defining ideal of the section of $D \to S'$ with the property that
${s'}^e=t$, where $e$ is the (local) degree of $f$. We then make an arbitrary 
lift
of $s'$ to an element $s$ of the defining ideal of the section of $D \to S$ (a
lifting now supposed to exist). We then can write $t$ as a power series
$\sum_{i=0}^\infty a_is^i$. As the deformation is supposed to be equisingular,
 we
have that $dt/ds$ is a unit times $s^{e-1}$. This forces $a_i=0$ for $0<i<e$, 
so that $t$ is of the form $a_0+s^eg$, where $g$ is a unit, as it is so in
$A'[[t]]$. Hence after possibly changing $s$ (but keeping its 
image in
$A'[[t]]$ the same) we may assume that $t=a_0+s^e$. However, the fact that
the reduced branch locus should be given by $t$ forces $a_0=0$, which gives
uniqueness. Existence is now trivial, as we may assume that $D_{|S'}$ has the
standard form given by $s^e=t$.

The case when $D$ is singular is very similar and left to the reader.
\qed
\enddemo

For any finite set $M$ of finite sets of integers $k_i\gs 2$,
$M=\{\{k_1,k_2,\dots,k_{n_1}\},\allowmathbreak\dots,\{\dots,k_{n_m}\}\}$,
such that
$2g-2=(2g(C)-2)d+\sum_i(k_i-1)$ we define $\sM_{g;M}(C/S)$ to be the
equisingular stratum of $\sM_{g;d}(C/S)$, with local degrees of ramification
$\gs 2$ over the points of the branch locus given by $M$.
Specializing even further, we let $\sM'$ be the open subset of
$\sM_{g;M}(\Pr^1/\Sp\cQ)$ for which $D$ is smooth, where $M$ has the form
$\{\{k_1,\dots,k_n\},\{2\},\{2\},\dots,\{2\}\}$. By Proposition~A.3,
$\sM'$ is smooth and so is then clearly also $\sM$, 
the closed substack for which
the point of the branch locus with ramification behavior
$\{k_1,\dots,k_n\}$
equals $\infty$. We now let $\sM^f$ be the stack over $\sM$ where the $D$ has
been provided with $n$ sections such that the fiberwise reduced ramification
locus over $\infty$ is the disjoint union of these sections together with a jet
of a local coordinate up to order $k_i$ at the $i$th section. The map
$\sM^f
\to \sM$ is evidently smooth and thus so is $\sM$.

Suppose we have the following data: A family $k\:D\to S$ of smooth proper 
connected
curves, $n$ disjoint sections $s_i$, a jet of order $k_i$ of a local
coordinate
at the $s_i$ and a polar part $p_i$ with respect to the coordinate of
order exactly
$k_i$ at the $s_i$. We then get a map $k_*\omega_{D/S} \to \sO_D$ given by
$\omega \mapsto \sum_i \Res_{s_i}p_i\omega$ and we assume that it is zero. This
condition may also be interpreted as follows. Consider the
line bundle $\sO\big(\sum_ik_is_i\big)$ on $D$. We have the natural
section
given by the fact that the divisor is effective, giving us a short exact 
sequence
$$
0\to{\sO_D}\to{\sO\bigg(\sum_ik_is_i\bigg)}\to{\sD}\to 0,
$$
where $\sD$ is the sheaf of polar parts of the appropriate orders at the $s_i$.
Pushing down by $k$ we get an exact sequence
$$
0 \to \sO_{S} \to k_*\sO\bigg(\sum_ik_is_i\bigg) \to k_*\sD \to
R^1k_*\sO_D.
$$
The $p_i$ give a section of $k_*\sD$, and the map 
$k_*\omega_{D/S} \to \sO_D$ is
the image of this section in $R^1k_*\sO_D$. The fact that it is zero thus gives
a vector bundle $\sE$ which is an extension of $\sO_S$ by $\sO_S$ mapping to
$k_*\sO\big(\sum_ik_is_i\big)$. As a further datum we require a splitting
of $\sE$
making it a trivial bundle. The two sections of $k_*\sO(\sum_ik_is_i)$ we 
get in this way have no common base points, and thus give an $S$-map 
$D \to \Pr^1_S$. We
further require that outside of $\infty$ this map only have simple 
singularities
(in the sense of $S$ being equal to the appropriate equisingular stratum).
We will let $\sN$ denote the moduli stack of such
data.

\proclaim{Proposition A.4}
The two stacks $\sM^f$ and $\sN$ are naturally isomorphic. In particular, $\sN$
is smooth.
\endproclaim

\demo{Proof}
By looking at the polar parts of the rational function $1/x$ on $\Pr^1$ 
we get a
map $\sM^f \to \sN$, and by looking at the map $D \to \Pr^1_{\sN}$ given by the
two sections of $k_*\sO\big(\sum_ik_is_i\big)$ we get a map for which all
of $\sN$ is
the  equisingular stratum given by
$\{\{k_1,\dots,k_n\},\{2\},\dots,\{2\}\}$, as outside of $\infty$ this is
built
into the definition of $\sN$, and at $\infty$ this follows from the 
specification
of the polar parts of $1/x$. This gives a map $\sN \to \sM^f$, and it is clear
that these two maps are each other's inverses.
\qed
\enddemo

Observe that the closure of $\sN$ represents the class
$[D_{\H}]$ defined in Sec.~5.3, while $\sM^f$ corresponds to 
$P\H_{g;k_1,\dots,k_n}$. Therefore, Proposition~A.4 implies that 
$[P\H_{g;k_1,\dots,k_n}]=[D_\H]$, as required in the proof of Lemma~5.5.


\Refs
\widestnumber\key{ACGH}

\ref \key{ACGH} 
\by E.~Arbarello, M.~Cornalba, P.~A.~Griffiths, and J.~Harris
\book Geometry of  Algebraic  Curves, vol.I
\publ Springer
\yr 1984
\endref

\ref\key{A1}
\by V.~I.~Arnold
\paper Topological classification
of complex trigonometric polynomials and combinatorics of graphs with an
identical number of vertices and edges
\jour Funct.~Anal. Appl.
\vol 30
\yr 1996
\pages 1--14
\finalinfo
\endref 

\ref\key{A2}
\by V.~I.~Arnold
\paper Critical points of functions and the classification of caustics
\jour Usp.~Math.~Nauk 
\vol 29 
\yr 1974 
\pages 243--244 
\finalinfo
\endref

\ref\key{CT}
\by M.~Crescimanno and W.~Taylor
\paper  Large $N$ phases of chiral QCD$_2$
\jour Nuclear Phys. B
\vol 437
\yr 1995 
\pages 3--24 
\finalinfo
\endref

\ref\key{D}
\by B.~A.~Dubrovin
\paper   Geometry of 2D topological field theories
\inbook Integrable Systems and Quantum Groups, Lecture Notes in Math.
\vol 1620
\yr 1996
\pages 120--348
\finalinfo
\endref

\ref\key{ELSV}
\by T.~Ekedahl, S.~K.~Lando, M.~Shapiro and A.~Vainshtein
\paper On Hurwitz numbers and Hodge integrals  
\jour C. R. Acad. Sci. Paris S\'er. I Math. 
\vol 328
\yr 1999
\pages 1175--1180
\finalinfo math.AG/9902104
\endref

\ref\key{FP1}
\by C.~Faber and R.~Pandharipande
\paper Hodge integrals and Gromov--Witten theory
\jour Invent. Math.
\vol 139
\yr 2000
\pages 173--199
\finalinfo math.AG/9810173
\endref

\ref\key{FP2}
\by C.~Faber and R.~Pandharipande
\paper Hodge integrals, partition matrices, and the $\lambda_g$-conjecture
\finalinfo math.AG/9908052
\endref

\ref\key{FnP}
\by C.~Fantechi and R.~Pandharipande
\paper Stable maps and branch divisors
\finalinfo math.AG/ 9905104
\endref

\ref \key{Fu} 
\by W.~Fulton
\book Intersection Theory
\publ Springer
\yr 1998
\endref

\ref\key{GL}
\by V.~Goryunov and S.~K.~Lando
\paper On enumeration of meromorphic functions on the line
\inbook The Arnoldfest, Fields Inst. Commun.
\vol 24
\yr 1999
\pages 209--223
\publ AMS
\publaddr Providence, RI
\endref

\ref\key{GJ1}
\by I.~P.~Goulden and D.~M.~Jackson
\paper Transitive factorization into transpositions
and holomorphic mappings on the sphere
\jour Proc. Amer. Math. Soc.
\vol 125
\yr 1997
\pages 51--60
\finalinfo 
\endref

\ref\key{GJ2}
\by I.~P.~Goulden and D.~M.~Jackson
\paper The combinatorial relationship between trees, cacti
and certain connection coefficients for the symmetric group
\jour Europ. J. Combin.
\vol 13
\yr 1992
\pages  357--365
\finalinfo 
\endref

\ref\key{GJ3}
\by I.~P.~Goulden and D.~M.~Jackson
\paper  A proof of a conjecture for the number of ramified
coverings of the sphere by the torus
\jour J. Combin. Theory Ser. A
\yr 1999
\vol 88
\pages 246--258
\finalinfo math.CO/9902009
\endref

\ref\key{GJVn}
\by I.~P.~Goulden, D.~M.~Jackson and A.~Vainshtein
\paper The number of ramified coverings of the sphere by the
torus and surfaces of higher genera
\jour Annals of Combin.
\vol 4
\yr 2000
\pages  27--46
\finalinfo math.CO/9902125
\endref

\ref\key{GJVk}
\by I.~P.~Goulden, D.~M.~Jackson and  R.~Vakil
\paper  The Gromow--Witten potential of a point,
Hurwitz numbers, and Hodge integrals 
\finalinfo math.AG/9910004
\endref

\ref\key{GP}
\by T.~Graber and R.~Pandharipande
\paper Localization of virtual classes
\jour Invent. Math.
\vol 135
\yr 1999
\pages  487--518
\finalinfo math.AG/9708001
\endref

\ref\key{GV}
\by T.~Graber and R.~Vakil
\paper   Hodge integrals and Hurwitz numbers via virtual localization
\finalinfo math.AG/0003028
\endref

\ref \key{HMo} 
\by J.~Harris and I.~Morrison
\book  Moduli of Curves
\publ Springer
\yr 1998
\endref

\ref\key{HMu}
\by J.~Harris and D.~Mumford
\paper On the {K}odaira dimension of the moduli space
  of curves
\jour Inv. Math.
\vol 67
\yr 1982
\pages  23--88
\finalinfo 
\endref

\ref\key{H}
\by A.~Hurwitz
\paper \"Uber Riemann'sche Fl\"achen
mit gegebenen Verzweigungspunkten
\jour  Math. Ann.
\vol 39
\yr 1891
\pages  1--60
\finalinfo 
\endref

\ref\key{K}
\by M.~Kontsevich
\paper Intersection theory
on the moduli space of curves and the matrix Airy function
\jour  Comm. Math. Phys.
\vol 147
\yr 1992
\pages  1--23
\finalinfo 
\endref

\ref\key{KM}
\by M.~Kontsevich and Yu.~I.~Manin
\paper Gromov--Witten
classes, quantum cohomology, and enumerative geometry
\jour  Comm. Math. Phys.
\vol 164
\yr 1994
\pages 525--562 
\finalinfo 
\endref

\ref\key{LZ}
\by S.~K.~Lando and D.~Zvonkine
\paper On multiplicities of the Lyashko--Looijenga mapping
on strata of the discriminant 
\jour  Funct. Anal. Appl
\vol 33
\yr 1999
\pages 178--188
\finalinfo 
\endref

\ref\key{L}
\by E.~Looijenga
\paper The complement of the bifurcation variety
of a simple singularity
\jour  Inv. Math.
\vol 23
\yr 1974
\pages 105--116
\finalinfo 
\endref

\ref\key{M1}
\by A.~D.~Mednykh
\paper Nonequivalent coverings of Riemann surfaces with a prescribed 
ramification type 
\jour  Sibirsk. Math. Zh.
\vol 25
\yr 1984
\pages 120--142
\finalinfo 
\endref

\ref\key{M2}
\by A.~D.~Mednykh
\paper  Branched coverings of Riemann
surfaces whose branch orders coincide with the multiplicity
\jour  Comm. in Algebra
\vol 18
\yr 1990
\pages 1517--1533
\finalinfo 
\endref

\ref \key{Mu} 
\by D.~Mumford
\book  Geometric Invariant Theory
\publ Springer
\yr 1965
\endref

\ref\key{N}
\by S.~Natanzon 
\paper  Topology of 2-dimensional
coverings and meromorphic functions on real
and complex algebraic curves
\jour  Selecta Mathematica Sovietica
\vol 12
\yr 1993
\pages 251--291
\finalinfo 
\endref

\ref\key{S}
\by V.~Strehl
\paper Minimal transitive products
of transpositions - the reconstruction of a proof by A.~Hurwitz 
\jour  Sem. Lothar. Combinat.
\vol 36
\yr 1996
\finalinfo Art. S37c
\endref

\ref\key{SSV}
\by B.~Shapiro, M.~Shapiro, and A.~Vainshtein
\paper  Ramified coverings of $S^2$ with one degenerate branching
point and enumeration of edge-oriented graphs
\publ AMS
\inbook Advances in Mathematical Sciences
\vol 34
\yr 1997
\pages 219--228
\publaddr Providence, RI
\endref

\ref\key{V1}
\by R.~Vakil
\paper Enumerative geometry of plane curves
of low genus 
\finalinfo math.AG/9803007
\endref

\ref\key{V2}
\by R.~Vakil
\paper Recursions, formulas,
and graphic-theoretic interpretations of ramified
coverings of the sphere by surfaces of
genus~$0$ and~$1$
\finalinfo math.AG/9812105
\endref

\ref\key{W}
\by E.~Witten
\paper  Two dimensional gravity and intersection
theory on moduli space
\inbook  Surveys in Differential Geometry
\publ Lehigh Univ.
\yr 1991
\pages 243--310
\finalinfo 
\endref

\ref\key{Z}
\by D.~Zvonkine
\paper Multiplicities of the Lyashko--Looijenga map
on its strata 
\jour  C.~R. Acad. Sci. Paris Ser.~I Math.
\vol 324
\yr 1997
\pages 1349--1353
\finalinfo 
\endref

\endRefs
\end{document}